\documentclass[twoside,11pt]{article}

\usepackage{amssymb}
\usepackage[enableskew]{youngtab}

\newcommand{\wt}{{\rm wt}}
\newcommand{\gr}{{\rm gr}}
\newcommand{\rk}{{\rm rk}}
\newcommand{\GKdim}{{\rm GKdim}}
\newcommand{\join}{\vee}
\newcommand{\meet}{\wedge}
\renewcommand{\k}{\Bbbk}

\renewcommand{\P}{{\mathbb P}}
\newcommand{\N}{{\mathbb N}}
\newcommand{\Z}{{\mathbb Z}}
\newcommand{\C}{{\mathbb C}}
\newcommand{\FF}{{\mathbb F}}
\newcommand{\A}{{\mathcal A}}
\newcommand{\T}{{\mathcal T}}
\newcommand{\M}{{\mathcal M}}
\renewcommand{\O}{{\mathcal O}}
\newcommand{\F}{{\mathcal F}}

\newcommand{\CC}{{\mathcal C}}
\renewcommand{\c}{{\mathbf c}}
\newcommand{\q}{{\mathbf q}}
\newcommand{\sm}{\underline{\rm sm}}
\newcommand{\la}{\langle}
\newcommand{\ra}{\rangle}

\newcommand{\inc}{\underline{{\rm inc}}}
\newcommand{\irr}{\underline{{\rm irr}}}
\newcommand{\opp}{{\rm opp}}
\newcommand{\tot}{{\rm tot}}
\newcommand{\tq}{\,|\,}
\newcommand{\proof}{\noindent {\it Proof.$\;$}}

\newcommand{\qed}{\hfill \rule{1.5mm}{1.5mm}}

\input xy
\xyoption{matrix}\xyoption{arrow}
\def\edge{\ar@{-}}
\def\dedge{\ar@{.}}

\newtheorem{theorem}{Theorem}[section]
\newtheorem{proposition}[theorem]{Proposition}
\newtheorem{definition}[theorem]{Definition}
\newtheorem{lemma}[theorem]{Lemma}
\newtheorem{corollary}[theorem]{Corollary}
\newtheorem{remark}[theorem]{Remark}
\newtheorem{example}[theorem]{Example}

\newtheorem{subtheorem}{Theorem}[subsection]
\newtheorem{subproposition}[subtheorem]{Proposition}
\newtheorem{subdefinition}[subtheorem]{Definition}
\newtheorem{sublemma}[subtheorem]{Lemma}

\newtheorem{subremark}[subtheorem]{Remark}

\newcommand{\titre}{Quantum analogues of 
Richardson varieties in the grassmannian and their toric degeneration.}

\begin{document}

\title{{\vspace{-1.5cm} \bf \titre}}
\author{L. RIGAL, P. ZADUNAISKY}
\date{}
\maketitle

\begin{abstract}
In the present paper, we are interested in natural quantum analogues of Richardson varieties in the type 
A grassmannians. To be more precise, the objects that we investigate are quantum analogues of the
homogeneous coordinate rings of Richardson varieties which appear naturally in 
the theory of quantum groups. 
Our point of view, here, is geometric: we are interested in the regularity properties of these 
{\em non-commutative varieties}, such as their irreducibility, normality, Cohen-Macaulayness... in the spirit 
of {\em non-commutative algebraic geometry}. A major step in our approach is to show that these algebras 
have the structure of an Algebra with a Straightening Law. From this, it follows that they degenerate 
to some quantum analogues of toric varieties.  \\
\end{abstract}

\section*{Introduction.}

\let\thefootnote\relax\footnotetext{2010 Mathematics Subject Classification. Primary: 16T20, 16S38, 17B37, 20G42. Secondary: 14M15.

Keywords and phrases: Quantum grassmannians, quantum Richardson varieties, quantum toric varieties, 
straightening laws, standard monomials, degeneration, Cohen-Macaulay, Gorenstein.}

Let ${\mathcal F}$ be a flag variety. As it 
is well known, the study of ${\mathcal F}$, both from the geometric 
and topological point of view, heavily relies on the study of its Schubert cells and Schubert varieties.
For example, the former give a stratification of ${\mathcal F}$ and the latter turn out to provide a nice 
understanding of the multiplicative structure of the cohomology ring of ${\mathcal F}$. In this context, 
it is 
important to understand how Schubert varieties and opposite Schubert varieties intersect. Such 
intersections 
are called Richardson varieties; they have been extensively studied for the last twenty years. 

It is beyond the scope of the present paper to give a complete overview of these studies and we will 
restrict 
ourselves to a short indicative list of works related to our own interests and considerations in the 
present 
work. Further references on the subject may be found by the interested reader in the few papers that we quote.

Early works on Richardson varieties include [D] and [R], where fundamental properties are studied, including 
their irreducibility. More recently, the extension of standard monomial theory to Richardson varieties was 
investigated in [LLit] in connection with some $K$-theoretic issues. Related results in the case of the type A 
grassmannians may be found in [KL]. 
Further, these varieties were used as central tools in [BL] in order to pursue the goal of providing 
a more geometric understanding of standard monomial theory. 

To finish this quick overview, let us mention the paper [M] where the existence of (semi-)toric degenerations 
of Richardson varieties is obtained, by means of representation theoretic methods based on canonical bases.\\

Beyond the classical case mentioned above, the theory of quantum groups provides natural analogues of Schubert 
varieties and, more generally, of Richardson varieties. Undoubtedly, the role of these 
quantum analogues in the study of quantum groups and quantum homogeneous spaces will be as central as it is 
in the classical setting. The objective of the present article is to establish some fundamental results 
regarding these objects,
in the case where the flag variety into consideration is a type A grassmannian.

To start with, let us briefly describe quantum Richardson varieties in the type A grassmannian case. 
Fix an arbitrary base field $\k$ and $q\in\k^\ast$. 
Let $m,n$ be integers such that $1 \le m \le n$. Consider the quantum analogue of the coordinate ring on 
the affine space of $n \times m$ matrices: $\O_q(M_{n,m}(\k))$ and let $\Pi_{m,n}\subseteq \N^m$ be the set of 
$m$-tuples 
$(i_1,\dots,i_m)$ of integers such that $1 \le i_1 < \dots < i_m \le n$ 
equipped with the obvious product order inhereted from $\N^m$.
To any element $I \in \Pi_{m,n}$, we may associate the quantum minor, denoted $[I]$, 
of  $\O_q(M_{n,m}(\k))$ built on the rows with index 
in $I$ and columns $1$ to $m$ of the generic matrix of $\O_q(M_{n,m}(\k))$. 
The subalgebra, $\O_q(G_{m,n}(\k))$, 
of $\O_q(M_{n,m}(\k))$ generated by these quantum minors is a natural analogue of the 
homogeneous  coordinate ring of the grassmannian with respect to the Pl\"ucker embedding. It is then 
easy to associate to any element $I\in\Pi_{m,n}$ a quantum Schubert and quantum opposite 
Schubert variety by considering the factor algebras 
$\O_q(G_{m,n}(\k))/\langle [K],\, K\in\Pi_{m,n},\, K\not\leq I\rangle$ and 
$\O_q(G_{m,n}(\k))/\langle [K],\, K\in\Pi_{m,n},\, K\not\geq I \rangle$, respectively. 
A natural analogue of the 
Richardson variety associated to a pair $(I,J)$ of elements of $\Pi_{m,n}$ being then defined as
the factor algebra 
$\O_q(G_{m,n}(\k))/\langle [K],\, K\in\Pi_{m,n},\, K \not\in [I,J]\rangle$. (Notice the abuse of language: these 
algebras are actually 
quantum analogues of homogeneous coordinate rings rather than quantum analogues of varieties.) 
For the convenience of the reader, we have included a short appendix at the end of the paper where 
classical Richardson varieties in type A grassmannians are described in some details. The material in this 
appendix also provides a justification of the above definition of quantum Richardson variety. \\

It is worth noting, at this stage, that quantum analogues of Schubert and Richardson varieties in partial flag 
varieties attached to simple Lie algebras of any type may be defined in a natural way (see
[LRes] for the case of Schubert varieties). 
This requires, however, a representation theoretic approach which forces to work under strong genericity 
hypotheses on the deformation parameter. The interested reader may consult the paper [Y1] by M. Yakimov 
for such an approach. 

In contrast, in the present paper, we adopt a method which allows us to work over any base field and 
with no assumption on the deformation parameter. It is inspired by the standard monomial theory which 
originates in works of De Concini, Eisenbud, Procesi, Lakshmibai, 
in the classical setting, and has its roots in early investigations by Hodge on grassmannians. 
The interested reader may consult, in particular, the following references: [DEP], [GL], [LRag]. 
Indeed, we first show that quantum analogues of Richardson varieties may be endowed with the structure 
of a (symmetric quantum) Algebra with a Straightening Law. That is, we exhibit a standard monomial basis 
for this algebra, built from a finite ordered set of generators and show that straightening and 
commutation relations satisfying certain very particular combinatorial constraints hold among these 
generators. We then deduce that quantum analogues of Richardson varieties may be filtered in
such a way that their associated graded ring is a tractable quantum analogue of coordinate ring of a 
toric variety. Hence proving that quantum analogues of Richardson varieties degenerate to  
quantum analogues of toric varieties.
We then derive from this important results of quantum analogues of Richardon varieties such as their 
irreducibility and Cohen-Macaulayness in the sense of Artin and Shelter. 

Actually, we work in a more general context: we first introduce the class of {\it symmetric quantum 
algebras with a straightening law} and then isolate a subclass, the objects of which can be degenerated
to quantum toric rings. We then show that this machinery applies to quantum Richardson 
varieties. Needless to say, since we 
are working in the non-commutative setting, the aforementioned commutation relations are a main issue 
in the present work. Actually, these relations are taken into account in the algebra with a staightening 
law structure.   \\

We finish this introduction by stressing two main remarks.

First, similar results have been obtained in [LR1] and [LR2] in the case of 
quantum analogues of Schubert varieties (which are special cases of quantum Richardson varieties). 
However, to deal with the case of Richardson varieties is significantly harder. For example, to prove 
that quantum analogues of Richardson varieties are integral domains already requires the full strength of 
the results in the present paper.  Namely, the fact that quantum analogues of Richardson varieties 
degenerate to quantum analogues of toric varieties. (A much more elementary proof of the fact that quantum
analogues of Schubert varities are integral domains is given in [LR2].)

Second, quantum Richardson varieties are investigated here from the point of view of non-commutative 
algebraic geometry. However, the results of [LR1] and [LR2] concerning the special case of quantum
Schubert varieties have been used in [LLR] to study the latter from a different point of view. 
Namely, the non-commutative structure on quantum flag varieties actually comes from a Poisson structure 
on the corresponding classical flag variety, in the spirit of deformation theory. Recent works have shed
light on the connection between certain torus-invariant prime ideals in quantum flag varieties, 
symplectic leaves of the corresponding classical objects and certain cell decompositions in their 
{\it totally positive} counterpart. The reader may also consult the references [GLL1], [GLL2], [Y2] and the survey article [LL]  for
more details on this point of view. We expect quantum Richardson varieties will be a strong tool to 
further study these connections. \\
 
The paper is organised as follows. In section \ref{distributive-lattices}, basic results that we need on 
distributive lattices are recalled. The material in this section is well known. We have 
summerized these results for the convenience of the reader. In section \ref{section-toric} a class
of quantum algebras associated with any distributive lattice is introduced. These are quotients of
quantum affine spaces by binomials defined on the basis of the attached lattice. They are natural quantum
analogues of coordinate rings of toric varieties. In section \ref{section-sqASL}, 
we introduce the notion of {\em Symetric quantum graded algebra with a straightening law}
(symetric quantum graded A.S.L., for short). It is a 
subclass of the class of quantum graded algebras with a straightening law defined in [LR1; Def. 1.1.1]
designed in order to show that the quantum Richardson varieties enjoy such an A.S.L. structure. A 
subclass of the class of symetric quantum graded A.S.L. is then defined,  in section 
\ref{section-condC-degeneration}, by means of an extra condition (C) which is imposed. 
This class includes the quantum toric varieties introduced in section 
\ref{section-toric} which are particularly simple examples. 
Actually much more is true: the quantum toric varieties are "essential" such 
examples in the sense that any algebra in this class degenerate to such an algebra. In section 
\ref{QARV} we reach our original motivation. We first show that quantum grassmannians are 
symetric quantum graded A.S.L. satisfying condition (C). It follows that the same holds for quantum 
Richardson varieties. Several important properties of the latter are then derived.

\paragraph{\em Notation and conventions.} 
Let $A$ be a ring, $\Pi$ an ordered set and $\iota \, : \, \Pi 
\longrightarrow A$
a map. To any finite increasing sequence $\alpha_1 \le \dots \le \alpha_t$, of length
$t\in\N^\ast$ of elements of $\Pi$ we associate the element $\iota(\alpha_1) \dots 
\iota(\alpha_t)$, 
which we call the standard monomial of $A$ associated to the sequence $\alpha_1 \le \dots \le 
\alpha_t$. We 
adopt the convention that there exists a unique increasing sequence of elements of 
$\Pi$ of length $0$ to
which we associate the standard monomial $1_A$. Hence, we have the familly of standard monomials on 
$\Pi$, 
that we will denote $\sm_\Pi(A)$. If any two distinct finite increasing sequences  of elements of 
$\Pi$ 
give rise to distinct standard monomials, there is a well defined notion of length for standard 
monomials. 
Namely, in this case, we define the length of a standard monomial as the length of its associated 
sequence.\\

Throughout, $\k$ denotes a field. If $S$ is a finite set, its cardinality will be  
denoted $|S|$.

\section{Reminder on distributive lattices.} \label{distributive-lattices}

For the convenience of the reader, we recall well known facts on distributive lattices that we will use 
all along the text. For this exposition, we essentially follow the paragraphs 3.1 to 3.4 of [S].
However, the interested reader may also refer to the foundational book [B] of G. Birkhoff. 

\paragraph{\em Ordered sets.} 
By an ordered set, we will always mean a set endowed with a partial ordering. 
Let $(\Pi,\le)$ be an ordered set. We denote by $\inc(\Pi\times\Pi)$ the subset of $\Pi\times\Pi$ 
of elements $(x,y)$ such that $x$ and $y$ are incomparable. 
The interval associated to a pair $(\alpha,\beta)$ of elements of $\Pi$ is defined by
$[\alpha,\beta]=\{\gamma\in\Pi \tq \alpha \le \gamma \le \beta\}$.
Clearly, it is non-empty if and only if $\alpha \le \beta$.
A subset $\Omega$ of $\Pi$ will be
called a $\Pi$-ideal (resp. $\Pi^\opp$-ideal) provided it satisfies the following condition: 
for all $\omega\in\Omega$ and all $\pi\in\Pi$, if $\pi\le\omega$, then $\pi\in\Omega$
(resp. for all $\omega\in\Omega$ and all $\pi\in\Pi$, if $\pi\ge\omega$, then $\pi\in\Omega$).
Let $(\Pi,\le)$ be a finite ordered set.
For any $x\in\Pi$, the rank of $x$, denoted $\rk(x)$, is defined to be the greatest 
integer $t$ such that there exists a strictly increasing sequence $x_0 < \dots < x_t=x$ in 
$\Pi$. Further, the rank of $\Pi$, denoted $\rk(\Pi)$, is defined by 
$\rk(\Pi)=\max\{\rk(x),\, x\in\Pi\}$.

\paragraph{\em Distributive lattices.}
A lattice is an ordered set 
$(\Pi,\le)$ satisfying the following condition: for any pair $(x,y)$ of elements of $\Pi$, there 
exists 
two elements $x \meet y$ and $x \join y$ in $\Pi$ such that $x \meet y \le x,y \le x \join y$, and 
for all $z \in \Pi$, if $z \le x,y$ (resp. $x,y \le
z$), $z \le x \meet y$ (resp. $x \join y \le z$); clearly, such elements are necessarily 
unique. 
Hence, if the ordered set $(\Pi,\le)$ is a lattice, we are given two maps
\[
\begin{array}{ccrcl}
\meet & : & \Pi \times \Pi & \longrightarrow & \Pi \cr
 & & (x,y) & \mapsto & x \meet y
\end{array}
\quad\mbox{and}\quad
\begin{array}{ccrcl}
\join & : & \Pi \times \Pi & \longrightarrow & \Pi \cr
 & & (x,y) & \mapsto & x \join y
\end{array}.
\]
A finite lattice is a lattice whose underlying set is finite. Clearly, a finite lattice has a unique 
minimal and a unique maximal element. 
The lattice $(\Pi,\le)$ is said to be distributive if it satisfies the following property:
for all $x,y,z \in \Pi$, $x \meet (y \join z)= (x \meet y) \join (x \meet z)$ or, equivalently, the 
property:
for all $x,y,z \in \Pi$, $x \join (y \meet z)= (x \join y) \meet (x \join z)$. 
A sub-lattice of $(\Pi,\le)$ is a subset $\Pi'$ of $\Pi$ endowed with the restriction of 
$\le$ which is stable under the maps $\meet$ and $\join$. A morphism of lattices is a morphism 
of ordered sets which commutes (in the obvious way) with the join and meet maps of each lattice. 

Let $\Pi$ be a lattice. 
An element $z\in\Pi$ is called join-irreducible provided it is not minimal 
and satisfies the following condition: if $x,y$ are 
elements of $\Pi$ such that $z = x \join y$, then either $z=x$, or 
$z=y$. We will denote by $\irr(\Pi)$ the set of join-irreducible elements of $\Pi$ and by 
$\irr^+(\Pi)$ the set of elements of $\Pi$ which are either join-irreducible or minimal. 
Then, we have the following celebrated structure theorem due to Birkhoff. The reader is refered to 
[B; Theo. 3, p.59] or [S; \S 3.4] for a proof of this statement.

\begin{theorem} -- {\bf (Birkhoff)} -- \label{birkhoff}
Let $\Pi$ be a finite distributive lattice and let $\Pi_0=\irr(\Pi)$. Then, 
$\Pi$ is isomorphic, as a distributive lattice, to $J(\Pi_0)$, where $J(\Pi_0)$ is the set
of $\Pi_0$-ideals of $\Pi_0$, ordered by inclusion. Further, the rank of $\Pi$ coincides with the 
cardinality of $\Pi_0$.
\end{theorem}

\paragraph{\em Finite chain products.}
Clearly, for all $d \in \N^\ast$, $\N^d$ endowed with the obvious product order, is a distributive 
lattice where, for
$i=(i_1,\dots,i_d)$ and $j=(j_1,\dots,j_d)$ in $\N^d$, $i \meet j =(\min\{i_1,j_1\},\dots,
\min\{i_d,j_d\})$ and   
$i \join j =(\max\{i_1,j_1\},\dots,\max\{i_d,j_d\})$. 

The following example will be of crucial importance. For any integer $p \ge 2$, we let 
$\CC_p=\{1,\dots,p\}$.
Now, consider $d \in \N^\ast$ and $d$ integers $n_1,\dots,n_d \ge 2$. The subset 
$\CC_{n_1} \times\dots \times \CC_{n_d}$ of 
$\N^d$ is clearly a (distributive) sub-lattice of $\N^d$. Such distributive lattices 
will be called finite chain products.
The following construction, associated to any finite chain product will be of central importance in 
the sequel. Fix $d \in\N^\ast$, integers $n_1,\dots,n_d \ge 2$ and put $N=2\max\{n_1,\dots,n_d\}+1$. 
We consider the map
\[
\begin{array}{ccrcl}
\omega & : &\CC_{n_1}\times\dots\times\CC_{n_d}& \longrightarrow & \N \cr
 & & (i_1,\dots,i_d) & \mapsto & \sum_{t=1}^d i_t N^{d-t}
\end{array}.
\]
Notice that, by definition, $\omega$ assigns to any element $I=(i_1,\dots,i_d) \in 
\CC_{n_1}\times\dots\times\CC_{n_d}$ 
the integer for which the coefficients of the $N$-adic expansion are the entries of $I$. 
Clearly, $\omega$ is strictly increasing and, in particular, injective. 

Lemma \ref{property-omega} records an important property of $\omega$ that we will need latter on. 
To state it, we need the following notation and terminology. 
Let $d$ be a positive integer and $I,J$ be two elements of $\N^d$. We denote by 
$I \sqcup J$ the union of $I$ and $J$ as multisets. That is, $I \sqcup J$ records the elements of $\N$ 
appearing as
coordinates in $I$ or $J$ together with the sum of their number of occurences in $I$ and their number of 
occurences in $J$. 

Further, an element $I=(i_1,\dots,i_d) \in \N^d$ will be termed increasing if $i_1 \le \dots \le i_d$.

\begin{lemma} -- \label{pre-property-omega}
Let $d$ be a positive integer and $s$ an integer such that $1\le s \le d$. Let 
$I=(i_t)_{1 \le t \le d},J=(j_t)_{1 \le t \le d},K=(k_t)_{1 \le t \le d},L=(l_t)_{1 \le t \le d}$
be increasing elements of $\N^ d$ such that $K \le I \le J \le L$ and $K \sqcup L = I \sqcup J$.  
The following holds: if $j_t=l_t$ for all $1 \le t <s$, then $i_t=k_t$ for all $1 \le t \le s$.
\end{lemma}

\proof The proof is by finite induction on $s$. When $s=1$, the hypothesis is empty and we must prove 
that $i_1=k_1$. Suppose, to the contrary, that $k_1<i_1$. Then, due to the fact that $I$ and $J$ are 
increasing and that 
$i_1 \le j_1$, no entry in $I$ and $J$ may equal $k_1$, contradicting $K \sqcup L = I \sqcup J$. Now, let 
$s$ be an 
integer such that $2 \le s \le d$ and suppose the result is true for all integers up to $s-1$. Suppose that $j_t=l_t$ 
for $1 \le t < s$.
Then, by the induction hypothesis, we have that $i_t=k_t$ for all $1 \le t \le s-1$. Suppose now that 
$k_s<i_s$.
Since $I$ and $J$ are increasing elements and $i_s \le j_s$, the only possible occurences of $k_s$ in $I$ 
(resp. $J$)
are among $i_1, \dots,i_{s-1}$, (resp. $j_1, \dots,j_{s-1}$). But, since 
$(k_1, \dots,k_{s-1})=(i_1, \dots,i_{s-1})$,  and $(l_1, \dots,l_{s-1})=(j_1, \dots,j_{s-1})$, this 
violates the identity
$K \sqcup L = I \sqcup J$. Hence, $k_s=i_s$ and the result holds for $s$.\qed

\begin{lemma} -- \label{property-omega}
Let $d$, $n_1,\dots,n_d$ and $\omega$ be as above.
Consider increasing elements $I,J,K,L$ of $\CC_{n_1}\times\dots\times\CC_{n_d}$ 
such that $K < I,J < L$ and $K \sqcup L = I \sqcup J$, then: \\
(i) $\omega(I)+\omega(J)\le \omega(K) + \omega(L)$;\\
(ii) $\omega(I)+\omega(J)=\omega(K) + \omega(L)$ iff $K= I \meet J$ and $L = I \join J$.
\end{lemma}

\proof We proceed in two steps.\\
1. Suppose, first, that  $K < I \le J < L$.
Since $J < L$, there is an integer $s$, $1 \le s \le d$, such that $j_t = l_t$ for $1 \le t < s$
and $j_s < l_s$. But then, by 
Lemma \ref{pre-property-omega}, we have $i_t=k_t$ for all $1 \le t \le s$. It follows 
that $i_t+j_t=k_t+l_t$ for $1 \le t <s$
and $i_s+j_s < k_s+l_s$.  On the other hand, by the choice of $N$, 
\[
\omega(I)+\omega(J)=\sum_{t=1}^d (i_t+j_t)N^{d-t}
\quad\mbox{and}\quad
\omega(K)+\omega(L)=\sum_{t=1}^d (k_t+l_t)N^{d-t}
\]
and these are the $N$-adic expansions of $\omega(I)+\omega(J)$ and $\omega(K)+\omega(L)$, respectively. 
It follows at once that 
$\omega(I)+\omega(J) < \omega(K)+\omega(L)$.\\
2. We get back to the hypothesis of the lemma's statement. 
We have $K \le I \meet J \le I , J \le I \join J \le L$. 
Since $K \sqcup L=I \sqcup J$, we have either 
$K < I \meet J \le I , J \le I \join J < L$, or $K = I  \meet J \le I , J \le I \join J = L$. 
In the first case, point 1 above  
shows that $\omega(I)+\omega(J) = \omega(I \meet J) + \omega(I \join J) < \omega(K)+\omega(L)$, 
while in the second case, we 
clearly have $\omega(I)+\omega(J) = \omega(I \meet J) + \omega(I \join J) = \omega(K)+\omega(L)$. 
This proves the claim.\qed\\ 

We conclude this section by the following definition, to be used latter in the text.

\begin{definition} -- \label{def-realisation}
Let $(\Pi,\le)$ be a finite ordered set which is a distributive lattice. A realisation 
of $(\Pi,\le)$ in a finite chain product is a datum  $(d;n_1,\dots,n_d;\iota)$ where $d\in \N^\ast$, 
$n_1,\dots,n_d\in\N\setminus\{0,1\}$ and $\iota$ is an injective morphism of lattices from $\Pi$ to 
$\CC_{n_1} \times\dots\times\CC_{n_d}$.
\end{definition}

We note, in passing, that any finite distributive lattice admits a realisation in a finite chain product.
This is an easy consequence of Theorem \ref{birkhoff}. Indeed, let $\Pi$ be a finite distributive 
lattice, and put $\Pi_0=\irr(\Pi)$. Then, Birkhoff's theorem yields an isomorphism of lattices 
$\Pi \cong J(\Pi_0)$ (in the above notation). On the other hand, the ordered set 
$({\mathcal P}(\Pi_0),\subseteq)$ of all subsets of $\Pi_0$ is a distributive lattice 
in which $J(\Pi_0)$ naturally embeds. 
It remains to notice that $({\mathcal P}(\Pi_0),\subseteq)$ and $\CC(2)^{|\Pi_0|}$ are isomorphic as 
distributive 
lattices. Indeed, composing all these lattice morphisms provides us with a lattice embedding of $\Pi$ in 
$\CC(2)^{|\Pi_0|}$.

\section{A class of quantum toric algebras.}\label{section-toric}

In the present section,
we introduce a class of algebras associated to the datum consisting of a distributive 
lattice $\Pi$ and two maps $\q\, : \, \Pi \times \Pi \longrightarrow \k^\ast \mbox{ and } 
\c\,:\, \inc(\Pi\times\Pi) \longrightarrow \k^\ast$. 

These algebras will turn out (see sections \ref{section-condC-degeneration} and \ref{QARV}) to be natural 
degenerations of quantum analogues of Richardson varieties. Further, we investigate their basic properties 
and, 
notably, we show that under some hypothesis they are integral domains. These results will turn out to be 
crucial 
later in the paper to derive properties of quantum Richardson varieties.\\

Let $(\Pi,\le)$ be a finite ordered set which is a distributive lattice. 
Suppose we are given maps 
\[
\begin{array}{ccrcl}
\q & : & \Pi \times \Pi & \longrightarrow & \k^\ast \cr
 & & (\alpha,\beta) & \mapsto & q_{\alpha,\beta}
\end{array}
\qquad\mbox{and}\qquad 
\begin{array}{ccrcl}
\c & : & \inc(\Pi\times\Pi) & \longrightarrow & \k^\ast \cr
 & & (\alpha,\beta) & \mapsto & c_{\alpha,\beta}
\end{array}.
\] 
To the data consisting of 
$(\Pi,\le)$, $\q$ and $\c$, we associate the $\k$-algebra, denoted $\A_{\Pi,\q,\c}$, with 
generators 
$X_\alpha$, $\alpha\in\Pi$, and relations: 
\[
X_\alpha X_\beta = q_{\alpha,\beta} X_\beta X_\alpha, \,\, 
\forall (\alpha,\beta) \in \Pi\times\Pi ,
\]
and 
\[
X_\alpha X_\beta = c_{\alpha,\beta} X_{\alpha\meet\beta} X_{\alpha\join\beta}, \,\, \forall (\alpha,\beta) \in
\inc(\Pi\times\Pi).
\]

It is clear that $\A_{\Pi,\q,\c}$ is endowed a $\N$-grading where canonical generators all have degree 
one. It follows that the map $\Pi \longrightarrow \A_{\Pi,\q,\c}$, $\alpha \mapsto X_\alpha$, is 
injective. 

\begin{remark} -- \label{toric-is-sqgrASL} \rm 
Let $\Pi$, $\q$ and $\c$ be as above. Further, assume that standard monomials on $\Pi$
are linearly independent elements of the $\k$-vector space $\A_{\Pi,\q,\c}$. 
Then, $\A_{\Pi,\q,\c}$, endowed with
its natural grading, is a quantum graded A.S.L. on 
$(\Pi,\le)$ in the sense of [LR1; Def. 1.1.1]. In particular,
standard monomials on $\Pi$ form a basis of the $\k$-vector space $\A_{\Pi,\q,\c}$ 
(see [LR1; Prop. 1.1.4]).
\end{remark}

\begin{remark} -- \label{relations-on-parameters-toric} \rm 
Let $\Pi$, $\q$ and $\c$ be as above. Further, assume that standard monomials on $\Pi$
are linearly independent elements of the $\k$-vector space $\A_{\Pi,\q,\c}$ (in particular, 
they must be non-zero).\\
(i) Let $(\alpha,\beta)\in\Pi\times\Pi$. Then, clearly, $X_\alpha X_\beta$ is the product of some 
standard monomial by a non-zero scalar. It follows that $X_\alpha X_\beta\neq 0$. \\ 
(ii) Let $\alpha\in\Pi$. We have the relation 
$X_\alpha X_\alpha=q_{\alpha\alpha}X_\alpha X_\alpha$, so that 
$q_{\alpha\alpha}=1$. Further, let $(\alpha,\beta)\in\Pi\times\Pi$, we have the relation
$X_\alpha X_\beta=q_{\alpha\beta}X_\beta X_\alpha=q_{\alpha\beta}q_{\beta\alpha} X_\alpha X_\beta$. 
It follows that $q_{\alpha\beta}q_{\beta\alpha}=1$. \\
(iii) Let $(\alpha,\beta)\in\inc(\Pi\times\Pi)$. We have the relation 
$X_\alpha X_\beta=c_{\alpha\beta} X_{\alpha\meet\beta}X_{\alpha\join\beta}$ and
$q_{\beta\alpha}X_\alpha X_\beta=X_\beta X_\alpha=c_{\beta\alpha} X_{\alpha\meet\beta} 
X_{\alpha\join\beta}$. 
Hence, arguing as above, we get $c_{\alpha\beta}=q_{\alpha\beta}c_{\beta\alpha}$.
\end{remark}

\begin{example} -- \label{commutative-case-toric} \rm 
In the above notation, we may consider the case where the maps $\q$ and $\c$ are constant, equal to 
$1$.
In this case, putting $\q={\bf 1}$ and $\c={\bf 1}$, we get the algebra $\A_{\Pi,{\bf 1},{\bf 1}}$ 
which is 
just the quotient of the 
(commutative) polynomial ring in the indeterminates $\{X_\alpha,\,\alpha\in\Pi\}$ by the ideal 
$\la X_\alpha X_\beta - X_{\alpha\meet\beta}X_{\alpha\join\beta} ,\, (\alpha,
\beta)\in\inc(\Pi\times\Pi) \ra$. 
This ring has 
been extensively studied. It has been proved that it is a (quantum) graded algebra with a 
straightening law. Hence  standard monomials on $\Pi$  form a $\k$-basis of 
$\A_{\Pi,{\bf 1},{\bf 1}}$. 
Further, 
$\A_{\Pi,{\bf 1},{\bf 1}}$ is an integral domain. All the relevant details may be found in Hibi's 
original 
paper; see
[H; p. 100]. 
\end{example}

Our next aim is to show that, under the hypothesis that standard monomials are linearly independent 
elements of $\A_{\Pi,\q,\c}$, then $\A_{\Pi,\q,\c}$ is an integral domain. For this, we need to 
introduce a kind of universal version of $\A_{\Pi,\q,\c}$, denoted $\A_\Pi$,  designed in such a 
way that (under convenient hypotheses), the algebras  
$\A_{\Pi,\q,\c}$ be quotients of $\A_\Pi$.\\

Let $\Pi$ be a finite ordered set. Consider the free $\k$-algebra, $F_\Pi$,  on the set 
\[
S_\Pi=\{X_\alpha,\,\alpha\in\Pi\} \cup
\{Q_{\alpha\beta},\, (\alpha,\beta)\in\Pi\times\Pi)\} \cup\{C_{\alpha\beta},\, (\alpha,
\beta)\in\inc(\Pi\times\Pi)\}. 
\]
There is an $\N$-grading on $F_\Pi$ for which elements of 
$\{X_\alpha,\,\alpha\in\Pi\}$ all have degree $1$ and elements of
$\{Q_{\alpha\beta},\, (\alpha,\beta)\in\Pi\times\Pi)\} \cup\{C_{\alpha\beta},\, (\alpha,
\beta)\in\inc(\Pi\times\Pi)\}$
all have degree zero. Now, consider the ideal $I_\Pi$ of $F_\Pi$ generated by the following elements:\\
(i) $Q_{\alpha\beta}Q_{\beta\alpha}-1$, $(\alpha,\beta)\in\Pi\times\Pi$;\\
(ii) $C_{\alpha\beta}=Q_{\alpha\beta}C_{\beta\alpha}$, $(\alpha,\beta)\in\inc(\Pi\times\Pi)$;\\
(iii) $Q_{\alpha\beta}a-aQ_{\alpha\beta}$, $(\alpha,\beta)\in\Pi\times\Pi$, $a \in S_\Pi$;\\
(iv) $C_{\alpha\beta}a-aC_{\alpha\beta}$, $(\alpha,\beta)\in\inc(\Pi\times\Pi)$, $a \in S_\Pi$;\\
(v) $X_\alpha X_\beta - Q_{\alpha\beta} X_\beta X_\alpha$, $(\alpha,\beta)\in\Pi\times\Pi$;\\
(vi) $X_\alpha X_\beta - C_{\alpha\beta} X_{\alpha\meet\beta}X_{\alpha\join\beta}$, $(\alpha,
\beta)\in\inc(\Pi\times\Pi)$.\\
We let
\[
\A_\Pi = F_\Pi/I_\Pi.
\] 
Clearly, $I_\Pi$ is generated by homogeneous elements of $F_\Pi$, so that $\A_\Pi$ 
inherits from $F_\Pi$ an $\N$-grading. Slightly abusing notation, we still denote by $X_\alpha$, 
$C_{\alpha\beta}$, ... the images of the corresponding elements of $F_\Pi$ under the canonical 
surjection 
$F_\Pi \longrightarrow \A_\Pi$.  Using the obvious map $\Pi \longrightarrow F_\Pi \longrightarrow 
\A_\Pi$, we get standard monomials on $\Pi$ in $\A_\Pi$. Recall that
\[
\sm_\Pi(\A_\Pi)
=\{1\} \cup \{X_{\alpha_1} \dots X_{\alpha_\ell}, \, \ell\in\N^\ast,\, \alpha_1 \le \dots \le 
\alpha_\ell \in \Pi\} .
\]

\begin{lemma} -- \label{specialisation-A-pi} Retain the above notation.\\
(i) Consider maps $\q\, : \, \Pi \times \Pi \longrightarrow \k^\ast$ and $\c\,:\, 
\inc(\Pi\times\Pi) \longrightarrow \k^\ast$. If the set of standard monomials on $\Pi$ is linearly 
independent in $\A_{\Pi,\q,\c}$, then there is a surjective $\k$-algebra 
morphism $\A_\Pi \longrightarrow \A_{\Pi,\q,\c}$ such that $X_\alpha \mapsto X_\alpha$, $C_\alpha 
\mapsto c_\alpha$ and $Q_{\alpha\beta} \mapsto q_{\alpha\beta}$. 
(This applies in particular when $\q={\bf 1}$ and $\c={\bf 1}$.)\\
(ii) The elements of $\sm_\Pi(\A_\Pi)$ are pairwise distincts.
\end{lemma}

\proof Point (i) follows at once from Remark \ref{relations-on-parameters-toric} and the universal 
property of free algebras. (See also Example \ref{commutative-case-toric}). 
Point (ii) follows from point (i) and the fact that standard monomials form a basis of 
$\A_{\Pi,{\bf 1},{\bf 1}}$. \qed\\

Point (ii) of Lemma \ref{specialisation-A-pi} shows that there is a well defined notion of length 
for standard  monomials on $\Pi$ in $\A_\Pi$ (see the end of the introduction). \\

We denote by $\M_\Pi$ the multiplicative sub-monoid of $\A_\Pi$ generated by the set 
$\{Q_{\alpha\beta},\, (\alpha,\beta)\in\Pi\times\Pi)\} \cup\{C_{\alpha\beta},\, (\alpha,
\beta)\in\inc(\Pi\times\Pi))\}$.

\begin{lemma} -- \label{lemma-straightening-in-A-pi}
Let $s\in\N^\ast$ and $\alpha_1,\dots,\alpha_s$ be $s$ elements in $\Pi$. Then, there 
exist an element $\xi\in\M_\Pi$ and a standard monomial $\underline{m}\in\sm_\Pi(\A_\Pi)$ 
of length $s$ such that 
$X_{\alpha_1} \dots X_{\alpha_s}=\xi\underline{m}$ in $\A_\Pi$.
\end{lemma}

\proof We need to recall the notion of depth of an element of the finite ordered set $\Pi$. Let 
$\alpha\in\Pi$. 
The depth of $\alpha$ is defined to be the greatest integer $t$ such that there exists a strickly 
increasing 
sequence $\alpha=\alpha_0 < \dots < \alpha_t$ in $\Pi$.

We prove the statement by induction on $s$. The result is trivial when $s=1$. Let us assume the 
result is 
true up to a certain integer $s$. We have to prove that any product of $s+1$ elements of 
$\{X_\alpha,\,\alpha\in\Pi\}$ equals a standard monomial of length $s+1$, up to multiplication by 
an element of $\M_\Pi$. To 
do this, we proceed by (finite) induction on the depth of the last of these $s+1$ elements.
Let $\alpha_1,\dots,\alpha_{s+1}$ be elements of $\Pi$ such that $\alpha_{s+1}$ has depth $0$. 
By the induction 
hypothesis, there exists $\xi\in\M_\Pi$ and $\beta_1 \le \dots \le \beta_s\in\Pi$ such that 
$X_{\alpha_1} \dots X_{\alpha_{s+1}}=\xi X_{\beta_1} \dots X_{\beta_s} X_{\alpha_{s+1}}$. On the 
other hand, $\alpha_{s+1}$ must be the unique maximal element of $\Pi$. Hence, 
$X_{\beta_1} \dots X_{\beta_s} X_{\alpha_{s+1}}$ is a standard monomial and we are done. Suppose 
now the result is true whenever the depth of the last of these $s+1$ elements does not exceed $p$, 
for some integer $p$.
Consider $\alpha_1,\dots,\alpha_{s+1}$, elements of $\Pi$ such that $\alpha_{s+1}$ has depth $p+1$.
By the first induction hypothesis, there exist $\xi\in\M_\Pi$ and $\beta_1 \le \dots \le 
\beta_s\in\Pi$ such that $X_{\alpha_1} \dots X_{\alpha_{s+1}}=\xi X_{\beta_1} \dots X_{\beta_s} 
X_{\alpha_{s+1}}$. Three cases may then occur. \\
1. If $\beta_s \le \alpha_{s+1}$, then $X_{\beta_1} \dots X_{\beta_s} X_{\alpha_{s+1}}$ is a 
standard monomial, and we are done.\\
2. If $\beta_s > \alpha_{s+1}$, then
$X_{\beta_s}X_{\alpha_{s+1}}=Q_{\beta_s\alpha_{s+1}}X_{\alpha_{s+1}}X_{\beta_s}$. Hence, 
$X_{\alpha_1} \dots X_{\alpha_{s+1}}=\xi X_{\beta_1} \dots X_{\beta_s} X_{\alpha_{s+1}}
=\xi Q_{\beta_s\alpha_{s+1}} X_{\beta_1} \dots X_{\beta_{s-1}} X_{\alpha_{s+1}}X_{\beta_s}$. But, 
since $\beta_s > \alpha_{s+1}$, $\beta_s$ has depth at most equal to $p$, so that we may apply the 
second induction hypothesis to the product $X_{\beta_1} \dots X_{\beta_{s-1}} 
X_{\alpha_{s+1}}X_{\beta_s}$. Again, we are done.\\
3. If $\beta_s$ and $\alpha_{s+1}$ are not comparable, then $X_{\beta_s}X_{\alpha_{s+1}}=
C_{\beta_s\alpha_{s+1}}X_{\alpha_{s+1}\meet\beta_s}X_{\alpha_{s+1}\join\beta_s}$. 
Using the fact that $\alpha_{s+1}\join\beta_s > \alpha_{s+1}$, the same argument as in the 
second case allows to conclude.

This finishes the proof.\qed

\begin{proposition} -- \label{regularity-toric}
Let $\gamma$ be an element of $\Pi$.\\
1. Let $\underline{m}$ be a standard monomial of $\A_\Pi$, of length $s\in\N$. Then, there exists a 
unique 
standard monomial $\underline{m}'$ of $\A_\Pi$ such that there exists $\xi\in\M_\Pi$ satisfying 
$X_\gamma\underline{m}=\xi\underline{m}'$.\\
2. The map $\phi_\gamma \, : \, \sm(\A_\Pi) \longrightarrow   
\sm(\A_\Pi)$, $\underline{m} \mapsto \underline{m}'$ is injective.\\
3. Consider maps $\q\, : \, \Pi \times \Pi \longrightarrow \k^\ast$ and $\c\,:\, \inc(\Pi\times\Pi) 
\longrightarrow \k^\ast$. If the set of standard monomials on $\Pi$ is linearly independent in 
$\A_{\Pi,\q,\c}$, then the element $X_\gamma$ of $\A_{\Pi,\q,\c}$ is regular.
\end{proposition}

\proof 1. The existence follows at once from Lemma \ref{lemma-straightening-in-A-pi}. Now, let 
$\underline{m}$ be a standard monomial of length $s\in\N$. Suppose there exist distinct 
standard monomials 
$\underline{m}'$ and $\underline{m}''$ of arbitrary length such that there exists $\xi',\xi'' 
\in\M_\Pi$ for which 
$X_\gamma\underline{m}=\xi'\underline{m}'=\xi''\underline{m}''$. The image of these relations under 
the 
projection $\A_\Pi \longrightarrow \A_{\Pi,{\bf 1},{\bf 1}}$ gives, 
in $\A_{\Pi,{\bf 1},{\bf 1}}$, an equality 
between two distinct standard monomials of $\A_{\Pi,{\bf 1},{\bf 1}}$. 
But, this contradicts the linear 
independence of standard monomials over $\Pi$ in $\A_{\Pi,{\bf 1},{\bf 1}}$ 
(see Example \ref{commutative-case-toric}). Hence, the required unicity.\\
2. Consider distinct standard monomials $\underline{m}_1$ and $\underline{m}_2$ of $\A_\Pi$. 
Suppose there exists 
a standard monomial $\underline{m}$ of $\A_\Pi$ and elements $\xi_1,\xi_2\in\M_\Pi$ such that 
$X_\gamma\underline{m}_1=\xi_1\underline{m}$ and $X_\gamma\underline{m}_2=\xi_2\underline{m}$. 
The image of these two relations under the projection $p_{{\bf 1},{\bf 1}} \, : \,
\A_\Pi \longrightarrow \A_{\Pi,{\bf 1},{\bf 1}}$ then provide an equality 
$X_\gamma p_{{\bf 1},{\bf 1}}(\underline{m}_1)=X_\gamma p_{{\bf 1},{\bf 1}}(\underline{m}_2)$. But,
clearly, $p_{{\bf 1},{\bf 1}}(\underline{m}_1)$ and $p_{{\bf 1},{\bf 1}}(\underline{m}_2)$ are two 
distinct 
standard monomials over $\Pi$ of $\A_{\Pi,{\bf 1},{\bf 1}}$. This contradicts the integrity of
$\A_{\Pi,{\bf 1},{\bf 1}}$ (see Example \ref{commutative-case-toric}). Hence, the map $\phi_\gamma$ 
is injective.\\
3. By Remark \ref{toric-is-sqgrASL}, the set of standard monomials on $\Pi$ forms a $\k$-basis 
of $\A_{\Pi,\q,\c}$. 
Now, using points 1 and 2 above and the canonical projection $\A_\Pi \longrightarrow  
\A_{\Pi,\q,\c}$, we see that left multiplication by $X_\gamma$ in $\A_{\Pi,\q,\c}$ is a map which, 
up to scalar multiplication by elements of $\k^\ast$, sends injectively the standard monomial 
$\k$-basis of $\A_{\Pi,\q,\c}$ into itself. It follows that left multiplication
by $X_\gamma$ in $\A_{\Pi,\q,\c}$ is injective; that is, $X_\gamma$ is left-regular. Since 
$X_\gamma$ commutes, up to multiplication by elements of $\k^\ast$ with any standard monomial, it 
follows that $X_\gamma$ is also right-regular. This completes the proof.\qed\\
 
Let us consider maps $\q\, : \, \Pi \times \Pi \longrightarrow \k^\ast$ and $\c\,:\, 
\inc(\Pi\times\Pi) \longrightarrow \k^\ast$ and assume that the set of standard monomials on $\Pi$ 
is linearly independent in 
$\A_{\Pi,\q,\c}$. By Proposition \ref{regularity-toric}, the elements $X_\gamma$, $\gamma\in\Pi$, 
of $\A_{\Pi,\q,\c}$ are 
regular and, clearly, they are normal. Hence, we may form the left quotient ring of 
$\A_{\Pi,\q,\c}$ 
with respect to the multiplicative set generated by $\{X_\gamma,\, \gamma\in\Pi\}$, that we denote 
by $\A_{\Pi,\q,\c}^\circ$, and we have a canonical injection
\[
\A_{\Pi,\q,\c} \longrightarrow \A_{\Pi,\q,\c}^\circ.
\] 
Hence, in order to prove that $\A_{\Pi,\q,\c}$ is an integral domain, it suffices to prove that 
$\A_{\Pi,\q,\c}^\circ$ is an 
integral domain. It turns out that we can do even better. In fact, $\A_{\Pi,\q,\c}^\circ$ is 
isomorphic to a quantum torus (that is, a quantum analogue of a Laurent polynomial ring), as we 
now proceed to show. \\

Denote 
by $\T_{\Pi,\q,\c}$ the $\k$-algebra generated by elements $X_\alpha^{\pm 1}$, 
$\alpha\in\irr^+(\Pi)$ subject to the relations $X_\alpha X_\beta=q_{\alpha\beta}X_\beta X_\alpha$, 
for all $\alpha,\beta\in\irr^ +(\Pi)$. 
Then, clearly,  $\T_{\Pi,\q,\c}$ is a quantum torus (see Remark \ref{relations-on-parameters-toric}.) 
In addition, there is a $\k$-algebra morphism as follows:
\[
\begin{array}{ccrcl}
j_{\Pi,\q,\c} &:& \T_{\Pi,\q,\c} & \longrightarrow & \A_{\Pi,\q,\c}^\circ \cr
 & & X_\gamma & \mapsto & X_\gamma
\end{array}.
\]

\begin{theorem} \label{isomorphism-torus-toric}
Let $\Pi$ be a finite ordered set which is a distributive lattice and consider $\q\, : \, \Pi 
\times \Pi \longrightarrow \k^\ast$ and $\c\,:\, \inc(\Pi\times\Pi) \longrightarrow \k^\ast$. 
Suppose, further, that the set of standard monomials on $\Pi$ is linearly independent in 
$\A_{\Pi,\q,\c}$. Then, $j_{\Pi,\q,\c}$ is
a $\k$-algebra isomorphism. In particular, $\A_{\Pi,\q,\c}$ is an integral domain.
\end{theorem}

\proof 
1. Surjectivity. Let $B$ be the $\k$-subalgebra of $\A_{\Pi,\q,\c}^\circ$ generated by the 
elements $X_\gamma^{\pm 1}$, $\gamma\in\irr^+(\Pi)$. To prove that $j_{\Pi,\q,\c}$ is surjective, 
it suffices to show that $B=\A_{\Pi,\q,\c}^\circ$. Let us show, by induction on the rank of 
$\gamma$, that for all
$\gamma\in\Pi$, $X_\gamma^{\pm 1}\in B$. If $\gamma$ has rank $0$, then $\gamma$ is the unique 
minimal element 
of $\Pi$, so that $\gamma\in\irr^+(\Pi)$. Hence $X_\gamma^{\pm 1}\in B$. 
Now, suppose the result true up to rank $p\in\N$ and consider $\gamma\in\Pi$, of rank $p+1$. 
If $\gamma$ is join-irreducible, then $\gamma\in\irr^+(\Pi)$, so that $X_\gamma^{\pm 1} \in B$.
Otherwise, there exist $\alpha,\beta\in\Pi$ such that $\gamma=\alpha\join\beta$ and 
$\alpha\meet\beta,\alpha,\beta < \gamma$.
Hence, the rank of the elements $\alpha\meet\beta,\alpha,\beta < \gamma$ is
less than or equal to $p$ and, $X_\alpha^{\pm 1},X_\beta^{\pm 1},X_{\alpha\meet\beta}^{\pm 1}
\in B$ by induction hypothesis. 
On the other hand, we have $X_\alpha X_\beta= c_{\alpha\beta} X_{\alpha\meet\beta}X_\gamma$
in $\A_{\Pi,\q,\c}$. 
Thus, $X_\gamma^{\pm 1}\in B$. This finishes the induction and the proof that $B=\A_{\Pi,\q,\c}$.\\
2. Injectivity. As is well known, $\T_{\Pi,\q,\c}$ is an integral domain of Gelfand-Kirillov 
dimension equal 
to the cardinality of $\irr^+(\Pi)$. On the other hand, as we saw earlier, $\A_{\Pi,\q,\c}$ is a 
quantum graded A.S.L. on $\Pi$ in the sense of [LR1], so that its Gelfand-Kirillov 
dimension 
equals the rank of $\Pi$ plus $1$ by [LR1; Prop. 1.1.5] (beware, in 
[LR1] a different convention was adopted for the rank). 
Suppose $j_{\Pi,\q,\c}$ has a non trivial kernel $J$, which hence 
contains a 
regular element, then by [KLen; Lemma 3.1, Prop. 3.15], we must have 
\[
\begin{array}{ccc}
\rk(\Pi)+1 & = & \GKdim(\A_{\Pi,\q,\c}) \cr
 & \le & \GKdim(\A_{\Pi,\q,\c}^\circ) \cr
 & = & \GKdim(\T_{\Pi,\q,\c}/J) \cr
 & < & \GKdim(\T_{\Pi,\q,\c}) \cr
 & = & |\irr^+(\Pi)| \cr 
 & = & |\irr(\Pi)|+1.
\end{array} 
\]
But, this contradicts Birkhoff's theorem (cf. Theorem \ref{birkhoff}). Hence $J=\{0\}$ and we are done. \\
As already mentioned, $\T_{\Pi,\q,\c}$ is an integral domain, so that $\A_{\Pi,\q,\c}$ is also an 
integral domain. \qed

\section{Symmetric quantum algebras with a straightening law.}\label{section-sqASL}

The aim of this section is to provide a general framework to study quantum Richardson varieties. More precisely, 
we introduce the notion of a symmetric quantum graded algebra with a straightening law on an ordered set $\Pi$. 
This is a subclass of the class of quantum graded algebras with a straghtening law introduced in [LR1]. The main 
point here is that certain specific quotients of symmetric quantum graded algebras with a straightening law 
inherit the same structure from the original algebra. Examples of algebras within this class include the quantum 
Richardson varieties as will be shown in section \ref{QARV}.

\begin{definition} -- \label{symetric-q-gr-ASL}
Let $A$ be an $\N$-graded $\k$-algebra, $(\Pi,\le)$ be a finite ordered set and $\Pi \longrightarrow A$
be a map whose image consists of homogeneous elements of $A$ of positive degree which generate $A$ as a 
$\k$-algebra. We say that $A$ is a symmetric quantum graded algebra with a straightening law on $\Pi$ 
(symmetric quantum A.S.L., for short)
if the following three conditions are satisfied:\\
(i) the set of standard monomials on $\Pi$ is a linearly independent set;\\
(ii) for any pair $(\alpha,\beta)$ of incomparable elements of $\Pi$, there exists a relation 
$\alpha\beta = \sum_{(\lambda,\mu)} c_{\lambda,\mu}^{\alpha,\beta}\lambda\mu$ where the sum extends 
over pairs 
$(\lambda,\mu)\in\Pi^2$, with $\lambda< \alpha,\beta < \mu$ and where, for such a pair 
$(\lambda,\mu)$, $c_{\lambda,\mu}^{\alpha,\beta}\in\k$; \\
(iii) for any pair $(\alpha,\beta)$ of elements of $\Pi$, there exists a relation 
$\alpha\beta - q_{\alpha\beta} \beta\alpha 
= \sum_{(\lambda,\mu)} d_{\lambda,\mu}^{\alpha,\beta}\lambda\mu$ where the sum extends over pairs 
$(\lambda,\mu)\in\Pi^2$, with $\lambda< \alpha,\beta < \mu$, where, for such a pair 
$(\lambda,\mu)$, $d_{\lambda,\mu}^{\alpha,\beta}\in\k$ and where $q_{\alpha\beta} \in\k^\ast$.
\end{definition}

The following remarks aim at clarifying this definition. 

\begin{remark} -- \label{Pi-imbeds} \rm 
We retain the notation of Definition \ref{symetric-q-gr-ASL}. \\
(i) By condition (i) of Definition \ref{symetric-q-gr-ASL}, the image of $\Pi$ under the map $\Pi 
\longrightarrow A$ must be linearly independent. It follows that the map $\Pi \longrightarrow A$ must be 
injective. For this reason, we will often identify $\Pi$ with its image in $A$.\\
(ii) Clearly, a symmetric quantum graded A.S.L. is a quantum graded algebra with 
a straightening law in the sense of [LR1; Def. 1.1.1]. In particular, standard monomials 
on $\Pi$ actually form a $\k$-basis of a symmetric quantum A.S.L. 
(see [LR1; Prop. 1.1.4]).
\end{remark}
 
Remark \ref{unique-or-not} clarifies the status of the relations required by conditions (ii) and (iii) 
above.

\begin{remark} -- \label{unique-or-not} \rm 
We retain the notation of Definition \ref{symetric-q-gr-ASL}.\\
1. {\em Straightening relations.} Given any pair $(\alpha,\beta)$ of incomparable elements of $\Pi$, by 
the linear independence of standard monomials, there is a unique relation as required 
in point (ii) of Definition \ref{symetric-q-gr-ASL}. 
It will be called the straightening relation associated to the pair $(\alpha,\beta)$.\\
2. {\em Commutation relations.} Let $(\alpha,\beta)$ be a pair of elements of $\Pi$. A relation as 
required in point
(iii) of Definition \ref{symetric-q-gr-ASL} need not be unique. Any such relation will be called a 
commutation
relation associated to the pair $(\alpha,\beta)$.
\end{remark}

The next remark shows that, for any symmetric
quantum graded A.S.L., the straightening and 
commutation relations provide, actually, a presentation of the algebra. It also states an easy 
consequence which will be used latter.

\begin{remark} -- \label{presentation} \rm
We retain the notation of Definition \ref{symetric-q-gr-ASL}.\\
1. Let $\k\la\Pi\ra$ be the free algebra on $\Pi$. Hence, $\k\la\Pi\ra$ is freely generated, as an 
algebra, by elements $X_\pi$, $\pi\in\Pi$. \\
1.1. Let $(\alpha,\beta)$ be a pair of incomparable elements of $\Pi$. By
hypothesis, there is a (unique) straightening relation  
$\alpha\beta = \sum_{(\lambda,\mu)} c_{\lambda,\mu}^{\alpha,\beta}\lambda\mu$ where the sum extends 
over pairs $(\lambda,\mu)\in\Pi^2$, with $\lambda< \alpha,\beta < \mu$ and 
$c_{\lambda,\mu}^{\alpha,\beta}\in\k$ are scalars. Put $S_{(\alpha,\beta)}=
X_\alpha X_\beta - \sum_{(\lambda,\mu)} c_{\lambda,\mu}^{\alpha,\beta}X_\lambda X_\mu \in\k\la\Pi\ra$.\\
1.2. Let $(\alpha,\beta)$ be any pair of elements of $\Pi$. By hypothesis, we can choose a (not necessarily unique) 
commutation relation $\alpha\beta - q_{\alpha\beta} \beta\alpha 
= \sum_{(\lambda,\mu)} d_{\lambda,\mu}^{\alpha,\beta}\lambda\mu$ where the sum extends over pairs 
$(\lambda,\mu)\in\Pi^2$, with $\lambda< \alpha,\beta < \mu$, where, for such a pair 
$(\lambda,\mu)$, $d_{\lambda,\mu}^{\alpha,\beta}\in\k$ and where $q_{\alpha\beta} \in\k^\ast$. For such a 
choice, let us put $C_{(\alpha,\beta)}=X_\alpha X_\beta - q_{\alpha\beta} X_\beta  X_\alpha 
- \sum_{(\lambda,\mu)} d_{\lambda,\mu}^{\alpha,\beta}X_\lambda X_\mu \in\k\la\Pi\ra$.\\
1.3. Let $I$ denote the ideal of $\k\la\Pi\ra$ generated by the elements $S_{(\alpha,\beta)}$ and
 $C_{(\alpha,\beta)}$ of 
points 1.1 and 1.2 above. Then, clearly, there is a surjective algebra morphism
$\k\la\Pi\ra /I \longrightarrow A$, $X_\pi \mapsto \pi$. There is also an obvious injective map 
$\Pi \longrightarrow \k\la\Pi\ra \longrightarrow \k\la\Pi\ra/I$ which allows to consider $\Pi$ as a 
subset of $\k\la\Pi\ra/I$. 
We want to show 
that (with its usual grading where canonical generators all have degree one), $\k\la\Pi\ra/I$, is
a quantum graded A.S.L. over $\Pi$. Notice first that the set of standard monomials on $\Pi$ 
in $\k\la\Pi\ra/I$ maps onto the set of standard monomials on $\Pi$ in $A$. Hence, the former has to be 
linearly independent, since the latter is. Second, observe that the straightening and commutation 
relations needed in $\k\la\Pi\ra/I$ indeed are available. It follows that $\k\la\Pi\ra/I$ is a symmetric 
quantum graded A.S.L. (with respect to its canonical grading). In addition, the set of standard monomials 
on $\Pi$ in $\k\la\Pi\ra/I$ form a basis and the projection  $\k\la\Pi\ra /I \longrightarrow A$, 
$X_\pi \mapsto \pi$ must be an isomorphism.\\
2. It follows from the algebra isomorphism $\k\la\Pi\ra /I \longrightarrow A$, $X_\pi \mapsto \pi$, that 
(appart from its original grading) $A$ can be endowed with an alternative grading where elements of 
$\Pi$ all are homogeneous of degree $1$ and that $A$ is a quantum graded A.S.L. over $\Pi$ with respect 
to this grading. From this, it follows that if $\pi_1,\dots,\pi_s$ are elements of
$\Pi$ ($s\in\N^ \ast$), the expression of the product $\pi_1 \dots \pi_s$ as a linear combination of
standard monomials involve only standard monomials which are ordered products of $s$ elements of $\Pi$.
\end{remark}

Next, we want to show that factors of a symmetric quantum graded A.S.L. by certain classes of ideals, 
arising from certain subsets of $\Pi$, inherit from $A$ a natural structure of symmetric quantum graded 
A.S.L. Recall the notion of $\Pi$-ideal and $\Pi^\opp$-ideal from section \ref{distributive-lattices}. \\

Recall from [LR1] that, 
given $\alpha\in\Pi$, we say that a standard monomial of $A$ on $\Pi$ involves $\alpha$ provided
it may be written as $\alpha_1\dots\alpha_s$ ($s\in\N^\ast$) with $\alpha_1 \le \dots \le \alpha_s \in \Pi$ 
and $\alpha\in\{\alpha_1,\dots,\alpha_s\}$. We then have the following result which will be useful 
latter. 

\begin{proposition} -- \label{basis-for-mixed-pi-ideals}
Let $A$ be an $\N$-graded $\k$-algebra, $\Pi$ a finite generating subset of the $\k$-algebra $A$ which 
consists in homogeneous elements of positive degree and suppose that $\Pi$ is ordered in such a way that 
$A$ be a symmetric quantum graded A.S.L. \\
1. Let $\Omega$ be a $\Pi$-ideal of $\Pi$. Then the set of standard monomials of $A$ on $\Pi$ involving 
an element of $\Omega$ is a $\k$-basis of the ideal $\la\Omega\ra$ of $A$.\\
2. Let $\Omega$ be a $\Pi^\opp$-ideal of $\Pi$. Then the set of standard monomials of $A$ on $\Pi$ 
involving an element of $\Omega$ is a $\k$-basis of the ideal $\la\Omega\ra$ of $A$.\\
3. Let $s\in\N^\ast$ and consider subsets $\Omega_1,\dots,\Omega_s$ of $\Pi$ which are either $\Pi$-ideals 
or $\Pi^\opp$-ideals of $\Pi$.  Then, the set of standard monomials of $A$ on $\Pi$ involving an 
element of $\cup_{1 \le i \le s} \Omega_i$ is a $\k$-basis of the ideal $\la\cup_{1 \le i \le s} 
\Omega_i\ra$ of $A$.
\end{proposition} 

\proof Recall that standard monomials of $A$ on $\Pi$ form a $\k$-basis of $A$.\\
1. This is Proposition 1.2.5 of [LR1].\\
2. This is proved in a similar fashion than point 1 above. Here is a sketch of proof, for the convenience 
of the reader.
First, arguing as in [LR1; Lemma 1.2.1], one can show that any element $\alpha\in\Pi$ is 
normal modulo 
the ideal of $A$ generated by $\{\pi\in\Pi,\, \pi >\alpha\}$. Second, if $<_\tot$ is any total order on 
$\Pi$ which respects 
$\le$ in the sense of [LR1; p. 676], then the elements of $\Pi$ increasingly ordered 
following 
$<_\tot$ form a normalising sequence, so that $\la\Omega\ra$ is the left ideal generated by $\Omega$.  
Then, the same 
proof as in [LR1; Prop. 1.2.5] leads to the result.\\
3. This follows from points 1 and 2 above.\qed

\begin{corollary} -- \label{sqgrASL-and-mixed-quotients}
Let $A$ be an $\N$-graded $\k$-algebra, $\Pi$ a finite generating subset of the $\k$-algebra $A$ which 
consists in 
homogeneous elements of positive degree 
and suppose that $\Pi$ is ordered in such a way that $A$ be a symmetric quantum graded 
A.S.L. Further, let $s\in\N^\ast$, $\Omega_1,\dots,\Omega_s\subseteq \Pi$ be either $\Pi$-ideals or 
$\Pi^\opp$-ideals 
of $\Pi$ and put $\Omega=\cup_{1 \le i \le s} \Omega_i$. Then, $A/\la\Omega\ra$ is a symmetric quantum 
graded A.S.L. 
on $\Pi\setminus\Omega$, equipped with the order induced from that of $\Pi$.
\end{corollary}

\proof By Proposition \ref{basis-for-mixed-pi-ideals}, it is clear that the obvious map
$\Pi\setminus\Omega \longrightarrow A \longrightarrow A/\la\Omega\ra$ is injective. In addition, $A/\la\Omega\ra$ 
clearly 
inherits an $\N$-grading from that of $A$ such that images of the above map are a generating set of the 
$\k$-algebra 
$A/\la\Omega\ra$  which are homogeneous of positive degree.
Further, using again Proposition \ref{basis-for-mixed-pi-ideals} and the fact that the set of standard 
monomials
of $A$ on $\Pi$ form a $\k$-basis of $A$, we get that standard monomials of $A/\la\Omega\ra$ on 
$\Pi\setminus\Omega$
form a $\k$-basis of $A/\la\Omega\ra$.
The existence of convenient straightening and commutation relations in $A/\la\Omega\ra$ follows
immediately from the existence of such relations in $A$.\qed

\section{Toric degeneration for a certain class of symmetric A.S.L.} \label{section-condC-degeneration}

In this section, we show that a certain class of quantum graded 
symmetric A.S.L. may be endowed with a filtration  
in such a way that the corresponding 
associated graded ring be a quantum toric algebra (in the sense of section 
\ref{section-toric} above). This class is determined by a combinatorial condition that the underlying 
ordered set has to satisfy. Hence, we will heavily rely on material of section 
\ref{distributive-lattices}. Roughly speaking, this class consists of symmetric quantum graded A.S.L.  
whose underlying poset may be realised (in the sense of Definition \ref{def-realisation}) by means of a 
sub-lattice of a chain product consisting of increasing elements of $\N^d$ and whose straightening and 
commutation relations must also satisfy certain combinatorial constraints. 

\begin{definition} --  \label{definition-condition-C}
Let $A$ be an $\N$-graded $\k$-algebra, $\Pi$ be a finite subset of $A$ which generates $A$ as a 
$\k$-algebra and consists of homogeneous elements of positive degree and $\le$ be an order on $\Pi$. 
Suppose $A$ is a 
symmetric quantum graded A.S.L. on $(\Pi,\le)$. We say that $A$ satisfies condition (C) provided the 
following 
holds:\\
(iv) $(\Pi,\le)$ is a distributive lattice and there exists $d\in\N^\ast$, integers $n_1,\dots,n_d \ge 2$ 
and a lattice
embedding $\iota \, : \, \Pi \longrightarrow \CC_{n_1}\times\dots\times\CC_{n_d}$ 
whose images are all increasing elements of $\N^d$;\\
(v) for any pair $(\alpha,\beta)$ of incomparable elements of $\Pi$, in the straightening relation 
associated 
to $(\alpha,\beta)$, we have $c_{\alpha\meet\beta,\alpha\join\beta}^{\alpha,\beta} \neq 0$ and, for all 
$(\lambda,\mu)$ such that $c_{\lambda,\mu}^{\alpha,\beta}\neq 0$, $\alpha \sqcup \beta = \lambda \sqcup 
\mu$ 
(indentifying elements of $\Pi$ with their image under $\iota$);\\
(vi) for any pair $(\alpha,\beta)$ of elements of $\Pi$, there exists a commutation relation associated 
to 
$(\alpha,\beta)$ such that $d_{\alpha\meet\beta,\alpha\join\beta}^{\alpha,\beta} = 0$ and, for all 
$(\lambda,\mu)$ such that $d_{\lambda,\mu}^{\alpha,\beta}\neq 0$, $\alpha \sqcup \beta = \lambda \sqcup 
\mu$ 
(indentifying elements of $\Pi$ with their image under $\iota$).
\end{definition}

\begin{remark} -- \rm \label{rem-after-def-C}
The present remark clarifies point (vi) in the conditions of Definition \ref{definition-condition-C}.\\
Assume $A$ is a symmetric quantum graded A.S.L. and retain the notation of Definition 
\ref{symetric-q-gr-ASL}.\\
1. Let $(\alpha,\beta)$ be a pair of comparable elements of $\Pi$. 
If $\alpha\beta - q_{\alpha\beta} \beta\alpha 
= \sum_{(\lambda,\mu)} d_{\lambda,\mu}^{\alpha,\beta}\lambda\mu$ is a commutation relation associated to 
$(\alpha,\beta)$, then the condition $d_{\alpha\meet\beta,\alpha\join\beta}^{\alpha,\beta} = 0$ is 
automatically 
satisfied.\\
2. Suppose $A$ satisfies condition (iv) and (v) of Definition \ref{definition-condition-C}. 
If $(\alpha,\beta)$ is a pair of incomparable elements of $\Pi$, then a commutation relation associated 
to 
this pair of the 
type required in Definition  \ref{definition-condition-C} always exists. Indeed, by condition (v), we 
have a 
straightening 
relation associated to both $(\alpha,\beta)$ and $(\beta,\alpha)$ with 
$c_{\alpha\meet\beta,\alpha\join\beta}^{\alpha,\beta} \neq 0$ and 
$c_{\alpha\meet\beta,\alpha\join\beta}^{\beta,\alpha} \neq 0$. Combining them in the obvious way, we get 
the desired expansion of $\alpha\beta-c_{\alpha\meet\beta,\alpha\join\beta}^{\alpha,\beta}
(c_{\alpha\meet\beta,\alpha\join\beta}^{\beta,\alpha})^{-1}\beta\alpha$.
\end{remark}

For the rest of this section, we fix the following notation. We let $A$ be an $\N$-graded $\k$-algebra, 
$\Pi$ be a subset of $A$ which generates $A$ as a $\k$-algebra and consists of homogeneous elements of 
positive 
degree and $\le$ be an order on $\Pi$ and we assume that $A$ is a 
symmetric quantum graded A.S.L. on $(\Pi,\le)$ which satisfies condition (C) for a fixed choice of a
realisation of $(\Pi,\le)$ by means of $d\in\N^\ast$, integers $n_1,\dots,n_d \ge 2$ and a lattice
embedding $\iota \, : \, \Pi \longrightarrow \CC_{n_1}\times\dots\times\CC_{n_d}$ 
whose images are all increasing elements of $\N^d$. We may then associate to the distributive lattice 
$\CC_{n_1}\times\dots\times\CC_{n_d}$ a map $\omega$ as defined in section \ref{distributive-lattices} 
and consider the composition
\[
\Pi \stackrel{\iota}{\longrightarrow} \CC_{n_1}\times\dots\times\CC_{n_d} 
\stackrel{\omega}{\longrightarrow} \N
\]
which we still denote $\omega$ for simplicity.  This map allows to associate to any element of $\Pi$, and 
more 
generally to 
any standard monomial, an integer which we call its weight. For this purpose, we let 
\[
M=\max\{\omega(\pi),\,\pi\in\Pi\}.
\]

\begin{definition} -- \label{definition-weight} Retain the above notation.\\
1. The weight of an element $\pi\in\Pi$ is the positive integer, denoted $\wt(\pi)$, and defined by:
\[
\wt(\pi)=M+1-\omega(\pi).
\]
2. Let $\pi_1 \le \dots \le \pi_s$ be elements of $\Pi$, the weight of the standard monomial 
$\pi_1 \dots \pi_s$ is the positive integer,
denoted $\wt(\pi_1 \dots \pi_s)$, and defined by $\wt(\pi_1 \dots \pi_s)=\wt(\pi_1)+\dots+\wt(\pi_s)$. 
Moreover, we let the 
weight of the standard monomial $1$ be $0$.
\end{definition}

The following lemma states some properties of weights.

\begin{lemma} -- \label{basic-property-weight} Retain the above notation.\\
(i) If $\alpha,\beta$ are elements of $\Pi$ such that $\alpha < \beta$, then 
$\wt(\alpha) > \wt(\beta)$.\\
(ii) Any element of $\Pi$ whose weight is one is maximal.\\
(iii) For elements $\alpha,\beta\in\Pi$, $\wt(\alpha)+\wt(\beta)=\wt(\alpha\meet\beta)+
\wt(\alpha\join\beta)$.\\
(iv) Let $(\alpha,\beta)$ be a pair of incomparable elements of $\Pi$. 
For any pair $(\lambda,\mu)$ different from $(\alpha\meet\beta,\alpha\join\beta)$ appearing in the 
straightening relation associated to $(\alpha,\beta)$, we have 
$\wt(\lambda\mu) < \wt(\alpha)+\wt(\beta)$.\\
(v) Let $(\alpha,\beta)$ be a pair of elements of $\Pi$ and consider any commutation relation associated 
to
this pair and satisfying condition (vi) in Definition \ref{definition-condition-C}. 
For any pair $(\lambda,\mu)$ appearing on the right hand side of this 
commutation relation, we have $\wt(\lambda\mu) < \wt(\alpha)+\wt(\beta)$.
\end{lemma}

\proof As noticed earlier, the map $\omega$ is strictly increasing, which proves (i) and (ii). 
Point (iii) is clear by point (ii) of Lemma \ref{property-omega}. Finally, 
 Lemma \ref{property-omega},  gives (iv) and (v). \qed\\

Using weights, 
we can  filter the $\k$-vector space $A$. For all $i\in\Z$, we denote by $\F_i$ the 
$\k$-subspace of $A$ with basis the set of standard monomials of weight less than or equal to $i$. 
Clearly, for $i < 0$, $\F_i=\{0\}$ and 
$\F_0=\k.1$. It is obvious that 
$\F=(\F_i)_{i\in\Z}$ is an ascending, exhaustive ({\it i.e.} $\cup_{i\in\Z} \F_i =A$) 
and separated ({\it i.e.} $\cap_{i\in\Z} \F_i =\{0\}$) filtration of the
$\k$-vector space $A$. The next proposition shows it is also a filtration of $A$ as a $\k$-algebra.

\begin{lemma} -- \label{pre-filtration}
Let $s\in\N^\ast$ and $\pi_1,\dots,\pi_s$ be elements of $\Pi$. 
Then, $\pi_1\dots\pi_s\in\F_{\sum_{1 \le i \le s} \wt(\pi_i)}$.
\end{lemma}

\proof We proceed by induction on $s$. The result is obvious if $s=1$. Fix now an integer $s\in\N^\ast$ 
and 
assume  that the result holds for $s$. We must show that, given any $\pi_1,\dots,\pi_{s+1}\in\Pi$,
the product $\pi_1\dots\pi_{s+1}$ belongs to $\F_{\sum_{1 \le i \le s+1} \wt(\pi_i)}$. For this, we 
proceed by  
induction on $\wt(\pi_{s+1})$. Suppose that $\wt(\pi_{s+1})=1$. Then $\pi_{s+1}$ is the unique maximal
element of $\Pi$ (recall that $\Pi$ is a distributive lattice). On the other hand, by the induction 
hypothesis
and point 2 of Remark \ref{presentation}, 
$\pi_1\dots\pi_s$  is a linear combination of standard monomials $\pi_1' \dots \pi_s '$
such that 
$\wt(\pi_1' \dots \pi_s ') = \wt(\pi_1' ) + \dots + \wt(\pi_s ') \le \wt(\pi_1) + \dots + \wt(\pi_s )$. 
Hence, 
$\pi_1 \dots \pi_{s+1}$ is a linear combination of terms $\pi_1' \dots \pi_s ' \pi_{s+1}$ which are all 
standard monomials of
weight $\wt(\pi_1' \dots \pi_s '\pi_{s+1})= \wt(\pi_1' ) + \dots + \wt(\pi_s ') +\wt(\pi_{s+1})
\le \wt(\pi_1) + \dots + \wt(\pi_s ) +\wt(\pi_{s+1})$ and we are done. Consider now an integer $r 
\in\N^\ast$, and assume 
the result is true for any product of $s+1$ elements of $\Pi$ whose last element has weight less than or 
equal to $r$. 
Consider $\pi_1, \dots, \pi_{s+1} \in \Pi$  such that $\wt(\pi_{s+1}) = r+1$. Again, by the first 
induction hypothesis and 
point 2 of Remark \ref{presentation}, $\pi_1 \dots \pi_s = \sum_j f_j \pi_{1,j} \dots \pi_{s,j}$,
where, for all $j$, $f_j\in\k$, $\pi_{1,j} \le \dots \le \pi_{s,j}\in\Pi$ and 
$\wt(\pi_{1,j} \dots \pi_{s,j}) = \wt(\pi_{1,j})+ \dots +\wt(\pi_{s,j}) \le \wt(\pi_1)+\dots+\wt(\pi_s)$.
Hence, we have 
\[
\pi_1 \dots \pi_s \pi_{s+1}= \sum_j f_j \pi_{1,j} \dots \pi_{s,j} \pi_{s+1}.
\]
The right hand side summands in the above equation fall into three possible cases.\\
{\em First case:} $\pi_{s,j} \le \pi_{s+1}$. For such a $j$, $\pi_{1,j} \dots \pi_{s,j} \pi_{s+1}$ is a standard monomial
of weight $\wt(\pi_{1,j}) + \dots + \wt(\pi_{s,j}) + \wt(\pi_{s+1}) \le \wt(\pi_1)+ \dots +\wt(\pi_s) + \wt(\pi_{s+1})$. So,  
$\pi_{1,j} \dots \pi_{s,j} \pi_{s+1} \in \F_{\sum_{1 \le i \le s+1} \wt(\pi_i)}$.\\
{\em Second case:} $\pi_{s,j} > \pi_{s+1}$. But, we have a commutation relation of the form 
\[
\pi_{s,j}\pi_{s+1}=d\pi_{s+1}\pi_{s,j}+\sum_{\lambda,\mu} d_{\lambda,\mu}\lambda\mu
\]
where the sum extends over pairs $(\lambda,\mu)$ of elements of $\Pi$ such that $\lambda< \pi_{s,j},\pi_{s+1} <\mu$.
Hence,
\[
\pi_{1,j} \dots \pi_{s,j} \pi_{s+1} 
= d\pi_{1,j} \dots \pi_{s-1,j}\pi_{s+1}\pi_{s,j}+\sum_{\lambda,\mu} d_{\lambda,\mu}\pi_{1,j} \dots \pi_{s-1,j}\lambda\mu.
\]
On the other hand, since $\pi_{s,j} , \mu > \pi_{s+1}$, we have $\wt(\pi_{s,j}) , \wt(\mu) < \wt(\pi_{s+1})$,
by Lemma \ref{basic-property-weight}. The second induction thus yields
\[
\pi_{1,j} \dots \pi_{s-1,j}\pi_{s+1}\pi_{s,j} \in \F_{\wt(\pi_{1,j})+ \dots +\wt(\pi_{s-1,j})+\wt(\pi_{s+1})+\wt(\pi_{s,j})}
\]
and 
\[
\pi_{1,j} \dots \pi_{s-1,j}\lambda\mu \in \F_{\wt(\pi_{1,j})+ \dots +\wt(\pi_{s-1,j})+\wt(\lambda)+\wt(\mu)}
\subseteq\F_{\wt(\pi_{1,j})+ \dots +\wt(\pi_{s-1,j})+\wt(\pi_{s+1})+\wt(\pi_{s,j})},
\]
the last inclusion being provided by Lemma \ref{basic-property-weight} which asserts that 
$\wt(\lambda)+\wt(\mu) \le \wt(\pi_{s,j})+\wt(\pi_{s+1})$. At this stage, we have that 
\[
\pi_{1,j} \dots \pi_{s,j} \pi_{s+1} \in \F_{\wt(\pi_{1,j})+ \dots +\wt(\pi_{s-1,j})+\wt(\pi_{s+1})+\wt(\pi_{s,j})}
\subseteq \F_{\wt(\pi_1)+ \dots +\wt(\pi_{s-1})+\wt(\pi_s)+\wt(\pi_{s+1})},
\]
since $\wt(\pi_{1,j})+ \dots +\wt(\pi_{s,j}) \le \wt(\pi_1)+ \dots +\wt(\pi_s)$.\\
{\em Third case:} $\pi_{s,j}$ and $\pi_{s+1}$ are not comparable. Proceeding as in case two by means of the 
straightening relation associated to the pair $(\pi_{s,j},\pi_{s+1})$, we get also that 
$\pi_{1,j} \dots \pi_{s,j} \pi_{s+1} \in \F_{\wt(\pi_1)+ \dots +\wt(\pi_{s-1})+\wt(\pi_s)+\wt(\pi_{s+1})}$.\\
Summing up the results of the three cases, we end up with  
$\pi_1 \dots \pi_{s+1} \in \F_{\sum_{1 \le i \le s+1} \wt(\pi_i)}$, as desired to complete the second induction and the proof.
\qed

\begin{proposition} -- \label{filtration} 
In the above notation, $\F$ is a filtration of the $\k$-algebra $A$.
\end{proposition}

\proof This follows at once from Lemma \ref{pre-filtration}.\qed\\

Our next aim in this section is to describe the associated graded ring, $\gr_\F(A)$, of $A$ with respect to the filtration $\F$. 
Hence, 
\[
\gr_\F(A) = \oplus_{i\in\Z} \F_i/\F_{i-1}=\oplus_{i\ge 0} \F_i/\F_{i-1}.
\]
Since $\F$ is separated and exhaustive, any non-zero element $x\in A$ has a {\em principal symbol} denoted $\gr(x)$. 
Namely, for any non-zero $x\in A$ there is a least integer $i\in\Z$ such that $x \in \F_i$; then we let
 $\gr(x)=x+\F_{i-1} \in \F_i/\F_{i-1}$.  Notice that, for all $x \in A\setminus\{0\}$, $\gr(x)\neq 0$.
 
It is clear, by definition of $\F$, that two distinct elements of $\Pi$ have distinct
principal symbols. Hence, $\gr(\Pi)=\{\gr(x),\,x\in\Pi\}\subseteq\gr_\F(A)$ identifies with $\Pi$ and inherits an order from it, together with a natural lattice structure such that $\Pi \longrightarrow \gr(\Pi)$, $x \mapsto \gr(x)$, be a lattice isomorphism. As a consequence, we get a natural realisation of $\gr(\Pi)$ in a finite chain product by means of that of $\Pi$.  \\

Our next aim is to show that $\gr_\F(A)$ is a symmetric quantum graded A.S.L. on $\gr(\Pi)$ satisfying condition (C). 
Notice that $\gr_\F(A)$ is naturally $\N$-graded. Notice, further, that for any standard monomial $x$ in $A$, then 
$\gr(x)$ is homogeneous of degree $\wt(x)$.

\begin{lemma} -- \label{pre-asgrad-ASL}
Retain the above notation. For $s\in\N^\ast$ and $\pi_1 \le \dots\le \pi_s \in \Pi$, we have
\[
\gr(\pi_1\dots\pi_s)=\gr(\pi_1)\dots\gr(\pi_s).
\]  
\end{lemma}

\proof If $w_1,\dots,w_s$ denote the weights of $\pi_1,\dots,\pi_s$, respectively, and if $w=w_1+\dots +w_s$, 
then by definition of the ring structure on $\gr(A)$, we have $\gr(\pi_1)\dots\gr(\pi_s)=\pi_1\dots\pi_s + \F_{w-1}$.
On the other hand, the standard monomial $\pi_1\dots\pi_s$ is in $\F_w\setminus \F_{w-1}$ since it has weight $w$;
so that $\gr(\pi_1\dots\pi_s)=\pi_1\dots\pi_s+\F_{w-1}$. Hence the result.\qed

\begin{proposition} -- \label{rel-in-agr} Retain the above notation.\\
1. The set $\gr(\Pi)$ generates $\gr_\F(A)$ as a $\k$-algebra.\\
2. Standard monomials of $\gr_\F(A)$ on $\gr(\Pi)$ are linearly independent.\\
3. For any pair $(\alpha,\beta)$ of incomparable elements of $\Pi$, there exists $c_{\alpha,\beta} 
\in\k^\ast$ such that
the following relation holds in $\gr_\F(A)$:
\[
\gr(\alpha)\gr(\beta)=c_{\alpha,\beta} \,(\gr(\alpha)\meet\gr(\beta))(\gr(\alpha)\join\gr(\beta)).
\]
4. For any pair $(\alpha,\beta)$ of elements of $\Pi$, there exists $q_{\alpha,\beta} \in\k^\ast$ such 
that the following relation holds in $\gr_\F(A)$:
\[
\gr(\alpha)\gr(\beta)=q_{\alpha,\beta} \,\gr(\beta)\gr(\alpha).
\]
5. The algebra $\gr_\F(A)$ is a symmetric quantum graded algebra on $\gr(\Pi)$ satisfying condition (C). 
\end{proposition}

\proof  Points 1 and 2 follow easily from Lemma \ref{pre-asgrad-ASL}.

Let $(\alpha,\beta)$ be a pair of incomparable elements of $\Pi$. By condition (v) of Definition 
\ref{definition-condition-C} and Lemma \ref{basic-property-weight}, we have that $\alpha\beta
\in \F_{\wt(\alpha)+\wt(\beta)}\setminus\F_{\wt(\alpha)+\wt(\beta)-1}$. So that, 
the straightening relation corresponding to these elements together with Lemma \ref{pre-asgrad-ASL} 
leads to a relation  
$\gr(\alpha)\gr(\beta)=\gr(\alpha\beta)=c_{\alpha,\beta}\gr((\alpha\meet\beta)(\alpha\join\beta))
=c_{\alpha,\beta}\gr(\alpha\meet\beta)\gr(\alpha\join\beta)
=c_{\alpha,\beta}(\gr(\alpha)\meet\gr(\beta))(\gr(\alpha)\join\gr(\beta))$, where $c_{\alpha,\beta} 
\in\k^\ast$. Further, by condition (vi) of Definition \ref{definition-condition-C}, and using the above,
we get point 4 of the proposition for such a pair. 

Let now $\beta \le \alpha \in \Pi$. Using condition (vi) of Definition 
\ref{definition-condition-C} for the pair $(\alpha,\beta)$, we see that $\alpha\beta
\in \F_{\wt(\alpha)+\wt(\beta)}\setminus\F_{\wt(\alpha)+\wt(\beta)-1}$. Arguing as above,
it follows that $\gr(\alpha)\gr(\beta)=\gr(\alpha\beta)=q_{\alpha,\beta}\gr(\beta\alpha)
=q_{\alpha,\beta}\gr(\beta)\gr(\alpha)$, where $q_{\alpha,\beta} 
\in\k^\ast$.

Points 3 and 4 are now proved. Point 5 is an obvious consequence.\qed\\

We now come to the description of $\gr_\F(A)$.

\begin{theorem} -- \label{theorem-degeneration-id}
Let $A$ be a symmetric quantum graded A.S.L. on $\Pi$ satisfying condition (C). 
Then, there exist maps $\q \, : \, \Pi \times \Pi \longrightarrow \k^\ast$, 
$\c \, : \, \inc(\Pi \times \Pi) \longrightarrow \k^\ast$ and a filtration $\F$ of the $\k$-algebra $A$ such that: \\
(i) $\gr_\F (A)$  
 is isomorphic (as a $\k$-algebra) to $\A_{\Pi,\q,\c}$; \\
(ii) standard monomials on $\Pi$ in $\A_{\Pi,\q,\c}$ are linearly 
independent.  \\
Further, $A$ is an integral domain.
\end{theorem}

\proof By the hypotheses, $A$ is an $\N$-graded $\k$-algebra, equipped with
a finite generating subset $\Pi$ consisting of homogeneous elements of positive degree, endowed 
with
an order such that $A$ is a symmetric quantum  graded A.S.L. Further, we suppose that $A$ satisfies 
condition (C). 
Then, by Proposition \ref{filtration}, $A$ admits a filtration $\F$ which gives rise to an 
associated graded ring $\gr_\F(A)$.
Now, using  Proposition \ref{rel-in-agr}  and the notation therein,
consider the maps $\q \, : \, \Pi \times \Pi \longrightarrow \k^\ast$,
$(\alpha,\beta) \mapsto q_{\alpha,\beta}$ and $\c \, : \, \inc(\Pi \times \Pi) \longrightarrow 
\k^\ast$,
$(\alpha,\beta) \mapsto c_{\alpha,\beta}$. By Proposition \ref{rel-in-agr} , there is a surjective 
$\k$-algebra morphism
\[
\begin{array}{ccrcl}
 & & \A_{\Pi,\q,\c}  & \longrightarrow & \gr_\F(A) \cr
 & & X_\alpha & \mapsto & \gr(\alpha)
\end{array} .
\] 
Further, Proposition \ref{rel-in-agr} ensures that standard monomials on $\gr(\Pi)$ are linearly 
independent elements of 
$\gr_\F(A)$, so that standard monomials on $\Pi$ are linearly independent elements of 
$\A_{\Pi,\q,\c}$. In addition, by 
Remark
\ref{toric-is-sqgrASL}, standard monomials on $\Pi$ form a basis of $\A_{\Pi,\q,\c}$. Hence, the above 
map is an 
isomorphism. From this and Theorem \ref{isomorphism-torus-toric}, we get that $\gr_\F(A)$ is an integral 
domain. By
well known results, it follows that $A$ is an integral domain.\qed\\

The following remark clarifies the non-unicity of commutation relations in symmetric 
quantum graded A.S.L satisfying condition (C).

\begin{remark} -- \rm  Let $A$ and $\Pi$ be as in Definition \ref{definition-condition-C} and assume $A$ 
satisfies 
condition (C) with respect to a realisation $(d;n_1,\dots,n_d;\iota)$ of $\Pi$ (see Definition 
\ref{def-realisation}). \\
1. Let $(\alpha,\beta)\in\Pi^ 2$. Suppose that in $A$ there exist two commutation relations as required 
by 
point (vi) of Definition \ref{definition-condition-C}, namely 
$\alpha\beta - q_{\alpha\beta} \beta\alpha 
= \sum_{(\lambda,\mu)} d_{\lambda,\mu}^{\alpha,\beta}\lambda\mu$ and 
$\alpha\beta - q_{\alpha\beta}' \beta\alpha 
= \sum_{(\lambda,\mu)} (d_{\lambda,\mu}^{\alpha,\beta})'\lambda\mu$. Then, by Proposition 
\ref{rel-in-agr} (and its proof), 
in $\gr_\F(A)$, we get the two 
relations: $\gr(\alpha)\gr(\beta) = q_{\alpha\beta} \gr(\beta)\gr(\alpha)$ and 
$\gr(\alpha)\gr(\beta) = q_{\alpha\beta}' \gr(\beta)\gr(\alpha)$. Hence, $0=(q_{\alpha\beta}-
q_{\alpha\beta}')\gr(\beta)\gr(\alpha)$. But, $\gr_\F(A)$ is an integral domain (see 
Proposition \ref{theorem-degeneration-id} and 
its proof). So, we must have $q_{\alpha\beta}-q_{\alpha\beta}'$. From this, we get that commutation 
relations as required by point (vi) of Definition
\ref{definition-condition-C} are actually unique.\\
2. By point 1 above, Remark \ref{unique-or-not} and Proposition \ref{rel-in-agr}, 
we may canonically associate to $A$ two maps:
 \[
 \begin{array}{ccrcl}
\q & : & \Pi\times\Pi & \longrightarrow & \k^\ast \cr
  & & (\alpha,\beta) & \mapsto & q_{\alpha\beta}
 \end{array}
\quad\mbox{and}\quad
 \begin{array}{ccrcl}
\c & : & \inc(\Pi\times\Pi) & \longrightarrow & \k^\ast \cr
 & & (\alpha,\beta) & \mapsto & c^{\alpha,\beta}_{\alpha\meet\beta,\alpha\join\beta}
 \end{array}
 \]
using the relevant non-zero scalars appearing in straightening relations and commutation relations as required by 
Definitions \ref{symetric-q-gr-ASL} and
\ref{definition-condition-C}, in such a way that the associated graded ring of $A$ with respect to 
the filtration $\F$ of Proposition \ref{filtration} be isomorphic to $\A_{\Pi,\q,\c}$. Using Remark 
\ref{relations-on-parameters-toric}, we then get that, \\
(i) $\forall\,(\alpha,\beta)\in\Pi\times\Pi$, $q_{\alpha\beta}q_{\beta\alpha}=1$ and 
$q_{\alpha\alpha}=1$;\\
(ii) $\forall\,(\alpha,\beta)\in\inc(\Pi\times\Pi$), 
$c^{\alpha,\beta}_{\alpha\meet\beta,\alpha\join\beta}
=q_{\alpha\beta}c^{\beta,\alpha}_{\alpha\meet\beta,\alpha\join\beta}$.
\end{remark}

We conclude this section by an easy observation, which will allow us latter to apply the above results to 
interesting 
classes of quantum algebras. \\

Let  $\Pi$ be an ordered set, $(\alpha,\beta)\in\Pi^2$ and $[\alpha,\beta]$ the corresponding interval 
(see section \ref{distributive-lattices}).
Letting
\[
\Pi_\alpha=\{\gamma\in\Pi \tq \alpha \not\le\gamma\}, 
\quad
\Pi^\beta=\{\gamma\in\Pi \tq \gamma \not\le\beta\} 
\quad\mbox{and}\quad
\Pi_\alpha^\beta=\Pi_\alpha\cup\Pi^\beta,
\]
we have that $[\alpha,\beta]=\Pi\setminus\Pi_\alpha^\beta$.

\begin{corollary} -- \label{condition-C-and-interval-ideals}
Let $A$ be a symmetric quantum graded A.S.L.  on the ordered set $\Pi$ and suppose that $A$ satisfies 
condition (C). 
Then, for any pair $(\alpha,\beta)$ of elements of $\Pi$ such that $\alpha\le\beta$, the quotient $\k$-algebra 
$A/\la\Pi_\alpha^\beta\ra$ is a symmetric quantum graded A.S.L.  on the ordered set $[\alpha,\beta]$, 
which satisfies
condition (C). In particular, $A/\la\Pi_\alpha^\beta\ra$ is an integral domain.
\end{corollary}

\proof Clearly, $\Pi_\alpha$ is a $\Pi$-ideal, $\Pi^\beta$ is a $\Pi^\opp$-ideal and the interval 
$[\alpha,\beta]$ is the complement of $\Pi_\alpha^\beta$ in $\Pi$. 
Hence, by Corollary \ref{sqgrASL-and-mixed-quotients} and its proof, 
the $\k$-algebra  $A/\la\Pi_\alpha^\beta\ra$ is a symmetric 
quantum graded A.S.L. on the ordered set $[\alpha,\beta]$ by means of the natural map
$[\alpha,\beta] \stackrel{\subseteq}{\longrightarrow}  A \stackrel{\rm can.proj.}{\longrightarrow} 
A/\la\Pi_\alpha^\beta\ra$. 
Further, $[\alpha,\beta]$ clearly inherits from $\Pi$ a distributive lattice structure as well as a
realisation as required by Definition \ref{definition-condition-C}.
It is then obvious that the 
straightening and commutation relations in $A/\la\Pi_\alpha^\beta\ra$, 
as required by Definition \ref{definition-condition-C},
may be obtained by applying the canonical projection 
$A \stackrel{\rm}{\longrightarrow} 
A/\la\Pi_\alpha^\beta\ra$ to the corresponding relations in $A$.
It remains to apply Theorem \ref{theorem-degeneration-id} to conclude that $A/\la\Pi_\alpha^\beta\ra$ 
is an integral domain.
\qed

\section{Quantum analogues of Richardson varieties.}\label{QARV}

In this section, we investigate quantum analogues of coordinate rings of Richardson varieties in the
grassmannians of type A. The final aim is to show that these are symmetric quantum graded A.S.L. 
satisfying condition (C) and to derive from this some of their important properties. \\

Consider positive integers $u,v$ and a scalar $q \in\k^\ast$.  
Following [LR1; sect. 3.1], we let $\O_q(M_{u,v}(\k))$ denote the quantum analogue of 
the affine coordinate ring of the space of $u \times v$ matrices with entries in $\k$. 
This is the $\k$-algebra with generators $X_{ij}$, $1 \le i \le u$ and $1 \le j \le v$ and relations as 
in [LR1; Def. 3.1.1]. If $u=v$, we put $\O_q(M_{v}(\k))=\O_q(M_{u,v}(\k))$.
Recall that there is a {\em transpose} automorphism of algebras ${\rm tr}_v \, : \, \O_q(M_{v}(\k)) \longrightarrow 
\O_q(M_{v}(\k))$, $X_{ij} \mapsto X_{ji}$. Recall in addition that, if $u',v'$ are positive integers 
such that $u'\le u$ and $v'\le v$, then the assignement $X_{ij} \mapsto X_{ij}$ defines an injective
algebra morphism from $\O_q(M_{u',v'}(\k))$ to $\O_q(M_{u,v}(\k))$. To any index sets $I,J$ 
of cardinality $t \le u,v$ with $I \subseteq \{1,\dots,u\}$ and $J \subseteq \{1,\dots,v\}$ we may 
associate a quantum minor, denoted $[I|J]$, and defined as in [LR1; section 3.1]. Then, it is well known 
that the transpose automorphism sends $[I|J]$ to $[J|I]$.  

Suppose now we are given integers $u, v$ such that $1 \le u \le v$. We let $\O_q(G_{u,v}(\k))$ denote the quantum analogue of the coordinate ring 
of the grassmannian of $u$-dimensional subspaces in $\k^v$. It is the subalgebra of $\O_q(M_{v,u}(\k))$ 
generated by the $u\times u$ quantum minors of $\O_q(M_{v,u}(\k))$. Notice that we adopt here a 
convention exchanging rows and columns with respect to the convention of [LR1] and [LR2]. However, 
embedding all the relevant algebras in $\O_q(M_{v}(\k))$ and using the transpose automorphism
introduced above shows that the two different conventions lead to isomorphic algebras. Hence, we are in 
position to use all of the results in the aforementioned papers.

Finally, we denote by $\Pi_{u,v}$
the subset of $\N^u$ of elements $(i_1,\dots,i_u)$ such that $1 \le i_1 < \dots < i_u \le v$ endowed
with the restriction of the natural product order of $\N^u$. 
It is easy to see that $\Pi_{u,v}$ is a (distributive) sub-lattice of $\N^u$. 
Clearly, an element $I=\{i_1 < \dots < i_u\}$ of $\Pi_{u,v}$ determines a
$u \times u$ quantum minor of $\O_q(M_{v,u}(\k))$ by sending $I$ to the minor
built on rows $i_1,\dots,i_u$; we denote this minor by $[I]$. 
The corresponding map $\Pi_{u,v} \longrightarrow \O_q(G_{u,v}(\k))$ turns
out to be injective (with image the canonical generators of $\O_q(G_{u,v}(\k))$).

\subsection{On the quantum grassmannians.}

We first prove some results on the quantum analogue of the coordinate ring of the grassmannians. 
Namely, we show that this $\k$-algebra is a symmetric quantum graded A.S.L. satisfying condition
(C). Part of this result is already contained in [LR1], where it is shown that this 
$\k$-algebra is a quantum graded A.S.L. in the sense of [LR1; Def. 1.1.1]. \\

Let $m \le n$ be positive integers and let $q\in\k^\ast$. Recall that to each elements $I \in\Pi_{m,n}$
we may associate the $m\times m$ quantum minor $[I]$ in 
$\O_q(G_{m,n}(\k))$ and that the map 
\[
\begin{array}{rcl} 
\Pi_{m,n} & \longrightarrow & \O_q(G_{m,n}(\k)) \cr
I & \mapsto & [I]
\end{array}
\]
is injective. We will often identify an element of $\Pi_{m,n}$ with its image in $\O_q(G_{m,n}(\k))$. 
By definition of $\O_q(G_{m,n}(\k))$, the set $\{[I],\, I\in\Pi_{m,n}\}$ is a set of generators of the 
$\k$-algebra $\O_q(G_{m,n}(\k))$ since any $m\times m$ quantum minor of $\O_q(M_{m,n}(\k))$ equals
$[I]$ for some $I\in\Pi_{m,n}$. Recall, further, that $\O_q(G_{m,n}(\k))$ has an $\N$-grading with 
respect
to which the elements $[I]$, $I\in\Pi_{m,n}$, are homogeneous of degree $1$.\\

It is proved in [LR1; sect. 3] that $\O_q(G_{m,n}(\k))$ is a quantum graded A.S.L.
on $\Pi_{m,n}$. More precisely, the following is proved:\\

\noindent
{\bf (ASL-1)} standard monomials on $\Pi_{m,n}$ form a basis of the $\k$-vector space 
$\O_q(G_{m,n}(\k))$; \\
{\bf (ASL-2)} for any $(I,J) \in\inc(\Pi_{m,n}\times\Pi_{m,n})$, there exists a (necessarily unique) 
relation 
\[
[I][J]=\sum_{(K,L)} c_{K,L}^{I,J} [K][L], 
\]
where the sum extends over pairs $(K,L)$ of elements of $\Pi_{m,n}$ such that $K \le L$ and $K < I,J$ and 
where, for such a pair, $c_{K,L}^{I,J}\in\k$.\\
{\bf (ASL-3)} for any $(I,J) \in\Pi_{m,n}\times\Pi_{m,n}$, there exists a relation 
\[
[I][J]-q^{f_{I,J}}[J][I]=\sum_{(K,L)} d_{K,L}^{I,J} [K][L], 
\]
where $f_{I,J}\in\Z$ and where the sum extends over pairs $(K,L)$ of elements of $\Pi_{m,n}$ such that 
$K \le L$ and $K < I,J$ and where, for such a pair, $d_{K,L}^{I,J}\in\k$.\\

We now establish a series of results allowing to get further information on the relations appearing in 
conditions (ASL-2) and (ASL-3) above. This will allow us to prove that $\O_q(G_{m,n}(\k))$ is a symmetric 
quantum
graded A.S.L. satisfying condition (C) on $\Pi_{m,n}$.
Notice that $\Pi_{m,n}$ is a distributive sub-lattice of $\N^m$ included in the finite chain product 
$\CC_n \times\dots\times \CC_n$ 
($m$ copies). Hence, the data $(m;n,\dots,n;\iota)$ is a realisation of $\Pi_{m,n}$ 
in a finite chain product in the 
sense of Definition \ref{def-realisation}, where 
$\iota \, : \, \Pi_{m,n} \stackrel{}{\longrightarrow} (\CC_n)^m$ is
the obvious inclusion map. Further, the elements of $\Pi_{m,n}$ are all increasing elements of $\N^m$. 
Hence, in order to prove that $\O_q(G_{m,n}(\k))$ satisfies condition (C), we will use this 
realisation.\\

In view of the defining relations of $\O_q(M_{n,m}(\k))$, it is easy to check 
that there exists a $\N^n$-grading on $\O_q(M_{n,m}(\k))$ such that, for $1 \le i \le n$ and 
$1 \le j \le m$,
$X_{ij}$ be of degree $\epsilon_i$. (Here, $\{\epsilon_1,\dots,\epsilon_n\}$ is the canonical 
basis of the $\Z$-module $\Z^n$.) Clearly, for 
$I=(i_1,\dots,i_m) \in \Pi_{m,n}$, 
$[I]$ is an homogeneous element of degree $\sum_{1 \le j \le m} \epsilon_{i_j}$. Hence, this $\N^n$-
grading on 
$\O_q(M_{n,m}(\k))$ induces by restriction an $\N^n$-grading on $\O_q(G_{m,n}(\k))$.

\begin{subproposition} -- \label{SR-CR-and-content}
Let $m\le n$ be positive integers and $q\in\k^ \ast$.\\
(i) Let $(I,J)\in\inc(\Pi_{m,n} \times \Pi_{m,n})$. In the straightening relation associated to $(I,J)$ 
in 
(ASL-2),
if $c_{K,L}^{I,J}\neq 0$, then $K \sqcup L = I \sqcup J$.\\
(ii) Let $(I,J)\in\Pi_{m,n} \times \Pi_{m,n}$. In any commutation relation associated to $(I,J)$ as in 
(ASL-3),
if $d_{K,L}^{I,J}\neq 0$, then $K \sqcup L = I \sqcup J$.
\end{subproposition}

\proof Using the $\N^n$-grading of $\O_q(G_{m,n}(\k))$ introduced above, the result follows from the  
linear independence of standard monomials in $\O_q(G_{m,n}(\k))$. \qed\\

We now need a technical lemma. \\

Notice that, for positive integers $u \le v$, there is an obvious one-to-one correspondence between 
elements of $\Pi_{u,v}$ and subsets of $\{1,\dots,v\}$ of cardinality $u$. Given $I \in\Pi_{u,v}$, we 
will call the corresponding subset the {\it underlying set of $I$}. 
  
\begin{sublemma} -- \label{muir-and-order-and-complement} 
Let $h,m,n$ be integers such that $1 \le h \le m \le n$.\\
(i) Let $I,J$ be elements of $\Pi_{m,2m}$. Denote by $I^c$ (resp. $J^c$) 
the element of $\Pi_{m,2m}$ whose underlying set
is $\{1,\dots,2m\}\setminus I$ (resp.  $\{1,\dots,2m\}\setminus J$). The following holds: if $I \le J$, 
then $I^c \ge J^c$. \\
(ii) Let $I,K\subseteq \Pi_{h,n}$ and $S \subseteq \{1,\dots,n\}$
of cardinality $m-h$ be such that $I \cap S=K \cap S=\emptyset$. 
The following holds: if $I \le K$ in $\Pi_{h,n}$, then $I \sqcup S \le K \sqcup S$ in $\Pi_{m,n}$.
\end{sublemma}

\proof Let $u\le v$ be positive integers. Observe that, for $I,J \in\Pi_{u,v}$, we have:
$I \le J$ if and only if, for all $1 \le i \le v$, $|I\cap\{1,\dots,i\}|\ge |J\cap\{1,\dots,i\}|$.
(Here, we identify an element in $\Pi_{u,v}$ and its underlying set.) The lemma follows at once.\qed

\begin{subproposition} -- \label{grassmannian-symetric-ASL}
Let $m\le n$ be positive integers and $q\in\k^ \ast$.\\
(i) Let $(I,J)\in\inc(\Pi_{m,n} \times \Pi_{m,n})$. 
In the straightening relation associated to $(I,J)$ in 
(ASL-2), if $c_{K,L}^{I,J}\neq 0$, then $K < I,J < L$.\\
(ii) Let $(I,J)\in\Pi_{m,n} \times \Pi_{m,n}$. There exists a commutation relation associated to $(I,J)$ 
as in (ASL-3) and such that, if $d_{K,L}^{I,J}\neq 0$, then $K < I,J < L$.
\end{subproposition}

\proof (i) To prove the first point, we proceed in several steps.\\
\noindent
{\em First case:} $n=2m$ and $I,J$ are elements of $\Pi_{m,2m}$ whose underlying sets have empty 
intersection. The straightening relation associated to $(I,J)$ in (ASL-2) is of the form
\[
[I][J]=\sum_{(K,L)} c_{K,L}^ {I,J} [K][L], 
\]
where the sum extends over pairs $(K,L)$ of elements of $\Pi_{m,2m}$ such that $K \le L$ and $K < I,J$ 
and 
where, for such a pair, $c_{K,L}^{I,J}\in\k$. 
Now, let $(K,L)$ be a pair such that $c_{K,L}^{I,J} \neq 0$. 
By point (i) of Proposition \ref{SR-CR-and-content}, 
the underlying sets of $K$ and $L$ are disjoint (and cover $\{1,\dots,2m\}$). 
Hence, in the notation of Lemma 
\ref{muir-and-order-and-complement}, $K^c=L$ and $L^c=K$.  
But, since $K < I$, Lemma \ref{muir-and-order-and-complement} gives that $L=K^c > I^c=J$. 
And we get that $L > I$ in the same way.
Thus point (i) holds in this case.\\
{\em Second case:} $n$ is arbitrary and $I,J$ are elements of $\Pi_{m,n}$ whose underlying sets have 
empty intersection 
(hence $n \ge 2m$). Letting $\A$ denote the subalgebra of $\O_q(M_{n,m}(\k))$ generated by the $X_{ij}$ 
with 
$i \in I \cup J$, we
have the obvious $\k$-algebra map $\O_q(M_{2m,m}(\k)) \longrightarrow \O_q(M_{n,m}(\k))$ wich is an 
embedding with image $\A$. Further, this map induces an embedding 
$\O_q(G_{m,2m}(\k)) \longrightarrow \O_q(G_{m,n}(\k))$. Now, a convenient choice of elements 
$I_0,J_0 \in \Pi_{m,2m}$ 
provides, using the first case above, a relation whose image under this last embedding is the desired 
relation in $\O_q(G_{m,n}(\k))$.\\
{\em Third case:} $n$ is arbitrary and $I,J$ are arbitrary elements of $\Pi_{m,n}$ ({\em i.e.} the 
general case). \\
Let $S$ stand for the intersection of the sets $I$ and $J$. So, there exists an integer $h$
such that $1 \le h \le m$ and elements $I_0,J_0\in\Pi_{h,n}$ satisfying $I=S \sqcup I_0$ and  
$J=S \sqcup J_0$. By definition, $I_0$ and $J_0$ do not intersect. Hence, the second case above provides 
a relation
\[
[I_0][J_0]=\sum c_{K_0,L_0} [K_0][L_0]
\]
in $\O_q(G_{h,n}(\k))$, where the sum extends on the pairs $(K_0,L_0)$ of elements of $\Pi_{h,n}$ such 
that $K_0 \le L_0$ and $K_0 < I_0,J_0 < L_0$ and where, for such a pair, $c_{K_0,L_0} \in \k$. 
Then, using 
Muir's law of extension of minors (see [LRu; Prop. 2.4]), we get a relation 
\[
[I][J]=[I_0\sqcup S][J_0\sqcup S]=\sum c_{K_0,L_0} [K_0\sqcup S][L_0\sqcup S]
\]
in $\O_q(G_{m,n}(\k))$. Using Lemma \ref{muir-and-order-and-complement}, we see that this is the 
desired relation. \\
The proof of point (i) is now complete.\\
(ii) The same reasonning as in point (i) may be applied.\qed

\begin{subproposition} -- \label{coeff-join-meet}
Let $m\le n$ be positive integers and $q\in\k^ \ast$.
For all $(I,J)\in\inc(\Pi_{m,n} \times \Pi_{m,n})$, there exists $e_{I,J}\in\Z$ such that, 
in the straightening relation associated to $(I,J)$ in 
(ASL-2), we have $c_{I\meet J,I \join J}^{I,J} = \pm q^{e_{I,J}}$.
\end{subproposition}

\proof We will have to use quantum grassmannians 
defined over arbitrary commutative integral domains as defined in [LR1; Def. 3.1.4] and
will make extensive use of [LR1; Rem. 3.1.5]. Put $\A=\Z[t^{\pm 1}]$.

Let $I,J$ be non comparable elements of $\Pi_{m,n}$. By [LR1; Theo. 3.3.8], in 
$\O_t(G_{m,n}(\A))$, we have a relation 
\[
[I][J]=\sum_{(K,L)} c_{K,L}(t) [K][L]
\]
where the sum extends over pairs $(K,L)$ of elements of $\Pi_{m,n}$ such that $K < I,J$ and $K \le L$ and 
where, for such a pair, $c_{K,L}(t)\in\A$. 

Let $u$ be any non-zero complex number and consider the obvious morphism 
$\Z[t^{\pm 1}] \longrightarrow \C[t^{\pm 1}] \longrightarrow \C$, where the second map is evaluation at 
$u$. It induces a map $\O_t(G_{m,n}(\Z[t^{\pm 1}])) \longrightarrow 
\O_t(G_{m,n}(\C[t^{\pm 1}])) \longrightarrow \O_u(G_{m,n}(\C))$
under which the image of the above straightening relation is 
\[
[I][J]=\sum_{(K,L)} c_{K,L}(u) [K][L] .
\]
Now, let $\delta= I \meet J$ and $\mu = I \join J$. 
We consider the quantum Schubert variety associated to $\delta$. This is the factor ring, 
$\O_u(G_{m,n}(\C))_\delta$ of $\O_u(G_{m,n}(\C))$
as defined in [LR1; Def. 3.1.7] and [LR2; Def. 1.1]. 
By Proposition \ref{SR-CR-and-content}, the image of the above relation in $\O_u(G_{m,n}(\C))_\delta$ is 
\[
[I][J]=c_{\delta,\mu}(u) [\delta][\mu] .
\]
But, by [LR2; Cor. 3.1.7], $\O_u(G_{m,n}(\C))_\delta$ is a domain and by 
[LR1; Cor. 3.4.5] it is a quantum graded A.S.L. on the subset 
$\{\pi\in\Pi_{m,n} \tq \pi \ge \delta\}$, so that 
the images of $[I]$ and $[J]$ in $\O_u(G_{m,n}(\C))_\delta$ are non-zero. 
Hence, we must have $c_{\delta,\mu}(u)\neq 0$. 

This shows that there exist integers $d,e$ such that $c_{\delta,\mu}(t)=dt^e$. 
Now, suppose that $d$ is divisible by some prime number $p\in\Z$, put $\FF_p=\Z/p\Z$ and consider 
the ring morphism $\Z[t^{\pm 1}] \longrightarrow \Z \longrightarrow \FF_p$, where the first map is 
evaluation at $1$ and the second is the canonical projection. We then 
have a natural ring morphism 
\[
\O_t(G_{m,n}(\Z[t^{\pm 1}])) \longrightarrow 
\O_1(G_{m,n}(\FF)) \longrightarrow 
\O_1(G_{m,n}(\FF))_\delta
\]
(the first arrow is obtained from the former morphism via [LR1; Rem. 3.1.5] and the 
second is the canonical projection). Applying this morphism 
to the above straightening relation would lead to the equantion $[I][J]=0$ in the (quantum) Schubert 
variety $\O_1(G_{m,n}(\FF))_\delta$ which would violate [LR2; Cor. 3.1.7] and 
[LR1; Cor. 3.4.5], as we argued above. So, $d=\pm 1$. 

Consider now the obvious ring morphism $\Z[t^ {\pm 1}] \longrightarrow \k$, $t \mapsto q$. 
It induces a map
$\O_t(G_{m,n}(\A)) \longrightarrow \O_q(G_{m,n}(\k))$ under which the image of the above 
straightening relation gives a straightening relation in $\O_q(G_{m,n}(\k))$ which establishes the claim.
\qed

\begin{subremark} -- \label{complement-coeff-join-meet} \rm 
Actually, Proposition \ref{coeff-join-meet} can be strengthened as follows. 
Let $m\le n$ be positive integers and $q\in\k^ \ast$.
For all $(I,J)\in\inc(\Pi_{m,n} \times \Pi_{m,n})$, there exists $e_{I,J}\in\Z$ such that, 
in the straightening relation associated to $(I,J)$ in 
(ASL-2), we have $c_{I\meet J,I \join J}^{I,J} = q^{e_{I,J}}$.
Indeed, it is known that in the classical case where $\k=\C$ and $q=1$, the coefficient 
$c_{I\meet J,I \join J}^{I,J}$ equals $1$; see [GL; Cor. 7.0.4, p. 236]. 
Using specialisation methods as in the proof of Proposition \ref{coeff-join-meet}, and using the notation 
therein, this easily leads to $d=1$. The result follows. 
\end{subremark}

\begin{subtheorem} -- \label{q-grass-is-sqASL-C}
Let $m \le n$ be positive integers and $q\in\k^\ast$. The $\k$-algebra $\O_q(G_{m,n}(\k))$
is a symmetric quantum graded A.S.L. on $\Pi_{m,n}$, which satisfies condition (C).
\end{subtheorem}

\proof We have already mentioned that standard monomials on $\Pi_{m,n}$ form a basis of the $\k$-vector 
space 
$\O_q(G_{m,n}(\k))$. Together with Proposition \ref{grassmannian-symetric-ASL}, this shows that 
$\O_q(G_{m,n}(\k))$ is a symmetric quantum graded A.S.L. Further, recall that we have the obvious 
realisation 
$(m;n,\dots,n;\iota)$ of $\Pi_{m,n}$ in a finite chain product discussed at the beginning of this 
subsection. Let us show that $\O_q(G_{m,n}(\k))$ satisfies condition (C) with respect to this 
realisation.
The existence of straightening relations as required by condition (v) of Definition 
\ref{definition-condition-C} is proved by Propositions \ref{SR-CR-and-content} and \ref{coeff-join-meet}. 
It remains to show, for each pair $(I,J)\in\Pi_{m,n}\times\Pi_{m,n}$, the existence of commutation 
relations 
as required by condition (vi) of \ref{definition-condition-C}. If $I$ and $J$ are not comparable, this 
follows from the above, by point 2 in Remark \ref{rem-after-def-C}. Suppose now $I$ and $J$ comparable. 
By Propositions \ref{SR-CR-and-content} and \ref{grassmannian-symetric-ASL}, 
there exists a commutation relation
as in (ASL-3) such that, if $d_{K,L}^{I,J}\neq 0$, then $K < I,J < L$ and $I \sqcup J = K \sqcup L$. 
But, since $I$ and $J$ are comparable, we have in particular that $d_{I\meet J,I \join J}^{I,J}= 0$. 
This finishes the proof.\qed

\subsection{Quantum Richardson varieties.}\label{ss-irreducibility}

In this subsection, we are interested in some quotients of the $\k$-algebra $\O_q(G_{m,n}(\k))$. These 
rings are natural quantum analogues of coordinate rings on Richardson varieties in the Grassmannian 
$G_{m,n}(\k)$ (see the introduction). \\

In the sequel, we  use the notation introduced before Corollary \ref{condition-C-and-interval-ideals}. 

\begin{subdefinition} -- Let $m \le n$ be positive integers and $q\in\k^\ast$. 
To each pair $(\alpha,\beta)$ of elements of $\Pi_{m,n}$ 
such that $\alpha\le\beta$, we associate the quantum analogue  of the homogeneous 
coordinate ring on the Richardson
variety corresponding to $(\alpha,\beta)$ also called, to simplify, the quantum Richardson variety 
associated to 
$(\alpha,\beta)$, defined as the quotient:
\[
\O_q(G_{m,n}(\k))/\la\Pi_\alpha^\beta\ra
\]   
of $\O_q(G_{m,n}(\k))$ by the ideal generated by the complement 
$\Pi_\alpha^\beta=\Pi_{m,n} \setminus [\alpha,\beta]$ of the
 interval $[\alpha,\beta]$ in $\Pi_{m,n}$.
\end{subdefinition}  

\begin{subtheorem} -- \label{qRV-are-sqgASL}
Let $m \le n$ be positive integers and $q\in\k^\ast$. For any pair $(\alpha,\beta)$ 
of elements of $\Pi_{m,n}$ such that $\alpha\le\beta$, the quantum Richardson variety 
$\O_q(G_{m,n}(\k))/\la\Pi_\alpha^\beta\ra$ is a symmetric quantum graded A.S.L. on the ordered set
$[\alpha,\beta]$, which satisfies condition (C). In particular, it is an integral domain.
\end{subtheorem}  

\proof This follows at once from Theorem \ref{q-grass-is-sqASL-C} and Corollary 
\ref{condition-C-and-interval-ideals}.\qed

\begin{subtheorem} -- \label{theo-GK}
Let $m \le n$ be positive integers and $q\in\k^\ast$. For any pair $(\alpha,\beta)$ 
of elements of $\Pi_{m,n}$ such that 
$\alpha=(\alpha_1,\dots,\alpha_m) \le (\beta_1,\dots,\beta_m)=\beta$, 
\[
\GKdim\, (\O_q(G_{m,n}(\k))/\la\Pi_\alpha^\beta\ra) = \sum_{k=1}^m (\beta_k-\alpha_k)+1.
\]
\end{subtheorem}  

\proof By Theorem \ref{qRV-are-sqgASL}, $\O_q(G_{m,n}(\k))/\la\Pi_\alpha^\beta\ra$ is
a symmetric quantum graded A.S.L. on the ordered set $[\alpha,\beta]$. In particular, it is a
quantum graded A.S.L. (in the sense of [LR1; Def.1.1.1]) on the ordered set $[\alpha,\beta]$. Thus, 
its Gelfand-Kirillov dimension is $r_{\alpha,\beta} + 1$, where $r_{\alpha,\beta}$ is the 
rank of the interval $[\alpha,\beta]$; see [LR1; Prop. 1.1.4]. (Beware, the rank of an ordered set as 
defined in [LR1] differ from the definition adopted in the present paper by one.)  
But, it is easy to show that $\rk([\alpha,\beta])=\sum_{k=1}^m (\beta_k-\alpha_k)$. We are done.\qed

\begin{subremark} -- \label{rem-sur-GKdim} \rm 
Let $m \le n$ be positive integers and $q\in\k^\ast$.
For any positive integer $d$, we let 
$S_d$ denote the $d$-th symmetric group. Identify $S_m \times S_{n-m}$ in the natural way with 
a subgroup of $S_n$. Then, it is well known that $\Pi_{m,n}$ identifies with the quotient
set $S_n/(S_m \times S_{n-m})$. More precisely, each coset of $S_n$ modulo $S_m \times S_{n-m}$ 
has a unique minimal element (with respect to the Bruhat ordering). Then, $\Pi_{m,n}$ identifies with the 
set of minimal coset representatives via the assignement sending $(\alpha_1,\dots,\alpha_m)$ to
the permutation $w_\alpha=(\alpha_1,\dots,\alpha_m,\beta_1,\dots,\beta_{n-m})$,
 where $\beta_1,\dots,\beta_{n-m}$
are the elements of the complement of $\{\alpha_1,\dots,\alpha_m\}$ in $\{1,\dots,n\}$ arranged in 
ascending order. It is then easily verified that, for each $\alpha\in\Pi_{m,n}$, the length of $w_\alpha$
is 
\[
\ell(w_\alpha)=\sum_{k=1}^m \alpha_k - \frac{1}{2}m(m+1).
\]
It follows that, for $\alpha\le\beta\in\Pi_{m,n}$, 
the formula of Theorem \ref{theo-GK} may be rewritten as:
\[
\GKdim\, (\O_q(G_{m,n}(\k))/\la\Pi_\alpha^\beta\ra) = \ell(w_\beta)-\ell(w_\alpha)+1.
\]
This applies in particular when the deformation parameter $q$ equals $1$ and shows that the 
Krull dimension of the homogeneous coordinate ring of the classical Richardson variety determined by 
$\alpha$ and $\beta$ (under the Pl\"ucker embedding) is $\ell(w_\beta)-\ell(w_\alpha)+1$ from which
it follows that, as a projective variety, it has dimension $\ell(w_\beta)-\ell(w_\alpha)$. Hence, we 
recover a well known result (see for example [LLit; Theo. 16]). 
\end{subremark}

We now investigate homological properties, namely the AS-Cohen-Macaulay and AS-Gorenstein properties.
For an overview on these notions as well as details useful in the sequel, the reader is refered to 
[LR1; section 2] where, moreover, a list of further references is available.

\begin{subtheorem} -- \label{theo-ASCM}
Let $m \le n$ be positive integers and $q\in\k^\ast$. For any pair $(\alpha,\beta)$ 
of elements of $\Pi_{m,n}$ such that 
$\alpha \le \beta$, the $\k$-algebra $\O_q(G_{m,n}(\k))/\la\Pi_\alpha^\beta\ra)$ is AS-Cohen-Macaulay.
\end{subtheorem}  

\proof By Theorem \ref{qRV-are-sqgASL}, $\O_q(G_{m,n}(\k))/\la\Pi_\alpha^\beta\ra$ is
a symmetric quantum graded A.S.L. on the ordered set $[\alpha,\beta]$. Further, the interval 
$[\alpha,\beta]$ is a distributive lattice, hence a wonderful poset in the sense of [LR1; Def. 2.2.3].
So, the result follows from a direct application of [LR1; Theo. 2.2.5].\qed 

\begin{subremark} -- \label{rem-ASG} \rm 
Let $m \le n$ be positive integers, $q\in\k^\ast$ and consider a pair $(\alpha,\beta)$ 
of elements of $\Pi_{m,n}$ such that $\alpha \le \beta$. 
We already argued that the $\N$-graded connected $\k$-algebra 
$\O_q(G_{m,n}(\k))/\la\Pi_\alpha^\beta\ra)$ is a 
quantum graded A.S.L in the sense of [LR1; Def.1.1.1]. It follows from Lemma 1.2.3 and Remark 2.1.4 of
[LR1]
that $\O_q(G_{m,n}(\k))/\la\Pi_\alpha^\beta\ra)$ is a noetherian algebra with enough normal elements (see 
[LR1; section 2.1]). Since, further, it is a domain and an AS-Cohen-Macaulay algebra, the fact that it is 
AS-Gorenstein or not can be read off from its Hilbert series. The interested reader may use section 4 of 
[LR1] to get details on this. However, the Hilbert series of $\O_q(G_{m,n}(\k))/\la\Pi_\alpha^\beta\ra)$
is independant of the value of $q\in\k^\ast$ since the standard monomials form a basis of this algebra 
consisting of homogeneous elements. Hence, the 
$\k$-algebra $\O_q(G_{m,n}(\k))/\la\Pi_\alpha^\beta\ra)$ is AS-Gorenstein if and only if the homogeneous 
coordinate ring of the Richardson variety $X_\alpha^ \beta$ is Gorenstein.
\end{subremark}

\begin{subremark} -- \label{rem-sur-ASCM} \rm 
It is worth noting, at this point, that in the proof of Theorem \ref{theo-ASCM} as well as in Remark 
\ref{rem-ASG}, the degeneration of 
quantum analogues of Richardson varieties to quantum toric varieties is not used. Indeed, the proof of 
these results only relies on the notion of a symmetric quantum graded A.S.L. as developed in section 
\ref{section-sqASL}. 
\end{subremark}

\begin{subremark} -- \label{} \rm 
We finish this work with a note concerning the normality of quantum Richardson varieties and, more generally, of 
symmetric quantum graded algebras with a straightening law satisfying condition (C).\\
1. Let $R$ be a noetherian domain, and denote by $Q$ its division ring of fractions. Then $R$ is a maximal order 
in $Q$ if, whenever $T$ is a subring of $Q$ such that $R \subseteq T \subseteq Q$ and there are elements 
$a, b \in R \setminus \{0\}$ with $aTb \subseteq R$, then $T = R$. Recall that if, in addition, $R$ is 
commutative, the above notion coincides with the usual notion of normality. Hence, we will say that a ring $R$ is 
normal if it is a noetherian domain which is a maximal order in its division ring of fractions.\\
2. Let $\Pi$ be a distributive lattice and consider maps 
$\q\, : \, \Pi \times \Pi \longrightarrow \k^\ast$ and 
$\c\,:\, \inc(\Pi\times\Pi) \longrightarrow \k^\ast$. Suppose further that standard monomials on $\Pi$ are 
linearly independent in $\A_{\Pi, \q, \c}$. As shown in Theorem \ref{isomorphism-torus-toric}, 
$\A_{\Pi, \q, \c}$ is an 
integral domain. Further, 
$\A_{\Pi, \q, \c}$ is noetherian by [LR1; Lemma 1.2.3], since it is a quantum graded A.S.L. It can be shown that 
$\A_{\Pi, \q, \c}$ is a normal ring. The natural way to prove this goes beyond the scope of the present paper. 
Actually $\A_{\Pi, \q, \c}$ turns out to belong to a class of natural quantum analogues of normal affine semigroup 
rings which can be shown to be normal domains. We intend to study this larger class of rings somewhere else.\\
3. Now, recall from Theorem \ref{theorem-degeneration-id} that if $A$ is a symmetric quantum graded A.S.L. satisfying 
condition (C), then $A$ can be filtered by an exhaustive separated filtration $\F$ such that $\gr_\F (A)$ is 
isomorphic to an algebra of type $\A_{\Pi, \q, \c}$ as in point 2 above. It follows that $A$ is a normal ring by 
standard results on the maximal order property concerning filtrations and associated graded rings (see [MR; 
Chapter X] or [McCR; \S 5.1.6]).\\
4. As a consequence of point 3 above, quantum Richardson varieties are normal rings (see Theorem 
\ref{qRV-are-sqgASL}).
\end{subremark}

\section*{Appendix.}

In the present appendix, we briefly recall the definition of Richardson varieties in the classical setting. 
For the sake of simplicity, we work over the field $\C$ of complex numbers and confine ourselves to 
Richardson varieties in the type A grassmannians. For this, we use results and notation from [GL; chap. 6].
Let $m,n$ be positive integers such that 
$1 \le m < n$. We let $G_{m,n}(\C)$ denote the grassmanian of $m$ dimensional subspaces of $\C^n$. 
Let $U \in G_{m,n}(\C)$. To an arbitrary basis $\{a_1,\dots,a_m\}$ of  $U$ we may associate the element 
$a_1 \wedge \dots \wedge a_m \in \bigwedge^m \C^n$ and its image $[a_1 \wedge \dots \wedge a_m]$ in
$\P(\bigwedge^m \C^n)$. Clearly, $[a_1 \wedge \dots \wedge a_m]$ is independent 
of the choice of the basis $\{a_1,\dots,a_m\}$ of $U$. Hence, we have defined a map
$p \, : \, G_{m,n}(\C)  \longrightarrow  \P(\bigwedge^m \C^n)$; 
the famous Pl\"ucker map. As is well known, $p$ is an embedding and its image is a closed subset of 
$\P(\bigwedge^m \C^n)$, so that $G_{m,n}(\C)$ acquires the structure of a projective variety. We put 
$G=SL(\C^n)$. Recall the natural action of  $G$ on $\P(\bigwedge^m \C^n)$. 

Let $\{e_1,\dots,e_n\}$ be the canonical basis of $\C^n$. Identifying $G$ with $SL_n(\C)$, we get the 
subgroups $T,B,B^-$ of $G$ corresponding to diagonal, upper-triangular and lower-triangular matrices.  
Denote by $\Pi_{m,n}$ the set of $m$-tuples $(i_1,\dots,i_m)$ of integers such that 
$1 \le i_1 < \dots < i_m \le n$. Further, for $I=(i_1,\dots,i_m) \in\Pi_{m,n}$, let 
$e_I =e_{i_1} \wedge \dots \wedge e_{i_m}$ (of course $[e_I]\in G_{m,n}$). 
It is easy to see that $\{[e_I], \,I\in\Pi_{m,n}\}$
is the set of fixed points of $\P(\bigwedge^m \C^n)$ for the natural action of $T$. 
For $I \in \Pi_{m,n}$, the Schubert cell associated to $I$ is defined as the $B$-orbit of $[e_I]$: 
$C_I=B.[e_I]$, while the Schubert variety, $X_I$, associated to $I$ is defined as the Zarisky closure of 
$C_I$ in $\P(\bigwedge^m \C^n)$: $X_I=\overline{C_I}$. It is not difficult to show that Schubert cells 
partition $G_{m,n}$. Of course, we can do the same looking at $B^-$-orbits rather than $B$-orbits. 
We then get opposite Schubert cells and varieties: $C^I=B^-.[e_I]$ and $X^I=\overline{C^I}$, for all
$I\in\Pi_{m,n}$. At this point, we may define the Richardson variety $X_J^I$ associated to a pair $(I,J)$ 
of elements of $\Pi_{m,n}$ as: $X_J^I=X^I \cap X_J$. 

Let $\C[p_I,\, I\in\Pi_{m,n}]$ be the homogeneous coordinate ring of $\P(\bigwedge^m \C^n)$. The image
by $p$ of the grassmannian $G_{m,n}(\C)$ is the set of points in $\P(\bigwedge^m \C^n)$ satisfying the 
Pl\"ucker relations. In addition, the algebra morphism from $\C[p_I,\, I\in\Pi_{m,n}]$ to
the polynomial ring  
$\C[X_{i,j},\, 1 \le i \le n, \, 1 \le j \le m]$ sending $p_I$ to the $m \times m$ 
(formal) minor $[I]$ of the generic matrix $(X_{ij})$ built on rows $i_1,\dots,i_m$, where
$I=(i_1,\dots,i_m)$, induces an isomorphism of algebras between the homogeneous coordinate ring of 
$p(G_{m,n}(\C))$ and the subalgebra of $\C[X_{i,j},\, 1 \le i \le n, \le 1 \le j \le m]$ generated by 
the elements $[I]$, $I \in\Pi_{m,n}$.

Let us now endow $\Pi_{m,n}$ with the obvious product order induced by $\Pi_{m,n} \subseteq \N^m$.
Consider $I,J\in\Pi_{m,n}$. 
It can be shown that the variety $X_J$
is the intersection of $p(G_{m,n}(\C))$ with the closed set of $\P(\bigwedge^m \C^n)$ defined by
the vanishing of $p_K$, for $K \not\leq J$. Similarly,
the variety $X^I$
is the intersection of $p(G_{m,n}(\C))$ with the closed set of $\P(\bigwedge^m \C^n)$ defined by
the vanishing of $p_K$, for $I \not\leq K$.
Hence, the variety $X_J^I$
is the intersection of $p(G_{m,n}(\C))$ with the closed set of $\P(\bigwedge^m \C^n)$ defined by
the vanishing of $p_K$, for $K\in\Pi_{m,n}\setminus [I,J]$. Notice that Schubert 
and opposite Schubert varieties are special cases of Richardson varieties since $\Pi_{m,n}$ has a lowest and a 
greatest element.

\section*{References.}

\noindent
[B]
G. Birkhoff. 
Lattice theory. Third edition. 
American Mathematical Society Colloquium Publications, Vol. 
XXV American Mathematical Society, Providence, R.I. 1967. 
\newline
[BL] M. Brion, V. Lakshmibai.
A geometric approach to standard monomial theory. 
Represent. Theory 7 (2003), 651--680.
\newline
[DEP] 
C. De Concini, D. Eisenbud, C. Procesi. 
Hodge algebras. Ast\'erisque, 91. Soci\'et\'e Math\'ematique de France, Paris, 1982.
\newline
[D] 
V. Deodhar. 
On some geometric aspects of Bruhat orderings. I. A finer decomposition of Bruhat cells.  
Invent. Math.  79  (1985),  no. 3, 499--511. 
\newline
[GLL1]
K.R. Goodearl, S. Launois, T.H. Lenagan. 
Totally nonnegative cells and Matrix Poisson varieties.
Adv. Math. 226 (2011) 779-826.
\newline
[GLL2]
K.R. Goodearl, S. Launois, T.H. Lenagan. 
Torus-invariant prime ideals in quantum matrices, totally nonnegative cells and symplectic leaves. 
To appear in Math. Z.
\newline
[GL] 
N. Gonciulea, V. Lakshmibai. 
Flag varieties. 
Travaux en Cours, 63. Hermann, 2001.
\newline
[H] T. Hibi. 
Distributive lattices, affine semigroup rings and algebras with straightening laws.  
Commutative algebra and combinatorics (Kyoto, 1985),  93--109, 
Adv. Stud. Pure Math., 11, North-Holland, Amsterdam, 1987.
\newline
[KLen]
G.R. Krause, T.H. Lenagan.
Growth of algebras and Gelfand-Kirillov dimension. 
Revised edition. Graduate Studies in Mathematics, 22. American Mathematical Society, Providence, 
RI, 2000.
\newline
[KL]
V. Kreiman, V. Lakshmibai. 
Richardson varieties in the Grassmannian.  
Contributions to automorphic forms, geometry, and number theory,  573–597, 
Johns Hopkins Univ. Press, Baltimore, MD, 2004. 
\newline
[LLit]
V. Lakshmibai, P. Littelmann. 
Richardson varieties and equivariant $K$-theory. 
Special issue celebrating the 80th birthday of Robert Steinberg.  
J. Algebra  260  (2003),  no. 1, 230--260. 
\newline
[LRag]
V. Lakshmibai, K. Raghavan. 
Standard monomial theory. Invariant theoretic approach. 
Encyclopaedia of Mathematical Sciences, 137. 
Invariant Theory and Algebraic Transformation Groups, 8. Springer-Verlag, Berlin, 2008.
\newline
[LRes]
V. Lakshmibai, N. Reshetikhin.
Quantum flag and Schubert schemes.  
Deformation theory and quantum groups with applications to mathematical physics (Amherst, MA, 1990), 
145--181. 
Contemp. Math., 134, Amer. Math. Soc., Providence, RI, 1992.
\newline
[LL]
S. Launois, T.H. Lenagan. 
From totally nonnegative matrices to quantum matrices and back, via Poisson geometry. 
To appear in the Proceedings of the Belfast Workshop on Algebra, Combinatorics and Dynamics 2009.
\newline
[LLR]
S. Launois, T.H. Lenagan, L. Rigal.
Prime ideals in the quantum Grassmannian.
Selecta Math. (N.S.) 13 (2008), no. 4, 697--725. 
\newline
[LR1]
T.H. Lenagan, L. Rigal. 
Quantum graded algebras with a straightening law and the AS-Cohen-Macaulay property for quantum 
determinantal rings and quantum Grassmannians.  J. Algebra  301  (2006),  no. 2, 670--702. 
\newline
[LR2]
T.H. Lenagan, L. Rigal.
Quantum analogues of Schubert varieties in the Grassmannian.  
Glasg. Math. J.  50  (2008),  no. 1, 55--70.
\newline
[LRu]
T.H. Lenagan, E.J. Russell. 
Cyclic orders on the quantum Grassmannian.  
Arab. J. Sci. Eng. Sect. C Theme Issues  33  (2008),  no. 2, 337--350.
\newline
[MR]
G. Maury, J. Raynaud. 
Ordres maximaux au sens de K. Asano. 
Lecture Notes in Mathematics, 808. 
Springer, Berlin, 1980.
\newline
[McCR] 
J.C. McConnell, J.C. Robson. 
Noncommutative Noetherian rings. 
Graduate Studies in Mathematics, 30. 
American Mathematical Society, Providence, RI, 2001.
\newline
[M] S. Morier-Genoud. 
Geometric lifting of the canonical basis and semitoric degenerations of Richardson varieties.  
Trans. Amer. Math. Soc.  360  (2008),  no. 1, 215--235.
\newline
[R]
R. Richardson.
Intersections of double cosets in algebraic groups.
Indag. Math. (N.S.) 3 (1992), no. 1, 69--77.
\newline
[S]
R.P. Stanley. Enumerative combinatorics. Vol. 1. 
Cambridge Studies in Advanced Mathematics, 49. Cambridge University Press, Cambridge, 1997. 
\newline
[Y1]
M. Yakimov.
A classification of $H$-primes of quantum partial flag varieties.  
Proc. Amer. Math. Soc.  138  (2010),  no. 4, 1249--1261.
\newline
[Y2]
M. Yakimov.
Invariant prime ideals in quantizations of nilpotent Lie algebras. 
Proc. Lond. Math. Soc. (3) 101 (2010), no. 2, 454--476.\\

\vskip .5cm

\noindent
Laurent RIGAL, \\
Universit\'e Paris 13, LAGA, UMR CNRS 7539, 99 avenue J.-B. Cl\'ement, 93430 Villetaneuse, 
France; e-mail: rigal@math.univ-paris13.fr\\

\noindent
Pablo ZADUNAISKY, \\
Universidad de Buenos Aires, FCEN, Departamento de Matem\'aticas, Ciudad Universitaria - Pabell\'on I - 
(C1428EGA) - Buenos Aires, Argentina; e-mail: pzadub@dm.uba.ar

\end{document}